\documentclass[11pt]{amsart}
\usepackage{comment}
\usepackage{amssymb}
\usepackage{graphicx}
\usepackage{mathtools}
\usepackage{subcaption}
\let\visiblecomments y
\usepackage{tikz-cd}
\usepackage{algorithm2e}

\newcommand{\pmap}{{\nrightarrow}}
\def\articletheorems{
\newtheorem{thm}{Theorem}[section]
\newtheorem{lem}[thm]{Lemma}

\newtheorem{defn}[thm]{Definition}

\newtheorem{prop}[thm]{Proposition}

}
\def\setof#1{\mbox{$\{\,#1\,\}$}}

\def\mathobj#1{\mbox{$#1$}}
\def\NN{\mathobj{\mathbb{N}}}
\def\RR{\mathobj{\mathbb{R}}}
\def\ZZ{\mathobj{\mathbb{Z}}}

\def\cA{\text{$\mathcal A$}}
\def\cC{\text{$\mathcal C$}}
\def\cK{\text{$\mathcal K$}}
\def\cP{\text{$\mathcal P$}}
\def\cQ{\text{$\mathcal Q$}}

\newcommand{\conv}{\protect\mbox{\rm conv\,}}

\newcommand{\cl}{\operatorname{cl}}
\newcommand{\bd}{\operatorname{bd}}
\newcommand{\Inv}{\operatorname{Inv}}
\newcommand{\Int}{\operatorname{int}}
\newcommand{\id}{\operatorname{id}}
\newcommand{\tr}{\operatorname{tr}}
\newcommand{\dom}{\operatorname{dom}}
\newcommand{\mto}{\multimap}
\newcommand{\perm}{\operatorname{Perm}}
\newcommand{\cycle}{\operatorname{Cycle}}
\newcommand{\Be}{\bar{B}_{\varepsilon}}

\articletheorems

\usepackage[cp1250]{inputenc}

\numberwithin{table}{section}

\begin{document}
\title{Conley index approach to sampled dynamics}
\author[B. Batko]{Bogdan Batko}
\address{Bogdan Batko,
         Division of Computational Mathematics,
         Faculty of Mathematics and Computer Science,
         Jagiellonian University,
         ul.~St. \L{}ojasiewicza 6, 30-348~Krak\'ow, Poland}
\email{bogdan.batko@uj.edu.pl}
\author[K. Mischaikow]{Konstantin Mischaikow}
\address{Konstantin Mischaikow,
         Department of Mathematics and BioMaPS Institute,
         Rutgers University,
         Piscataway, NJ 08854, USA}
\email{mischaik@math.rutgers.edu}
\author[M. Mrozek]{Marian Mrozek}
\address{Marian Mrozek,
         Division of Computational Mathematics,
         Faculty of Mathematics and Computer Science,
         Jagiellonian University,
         ul.~St. \L{}ojasiewicza 6, 30-348~Krak\'ow, Poland}
\email{marian.mrozek@uj.edu.pl}
\author[{M. Przybylski}]{Mateusz Przybylski}
\address{Mateusz Przybylski,
         Division of Computational Mathematics,
         Faculty of Mathematics and Computer Science,
         Jagiellonian University,
         ul.~St. \L{}ojasiewicza 6, 30-348~Krak\'ow, Poland}
\email{Mateusz.Przybylski@im.uj.edu.pl}
\date{today}
\date{\today}
\thanks{
   BB, MM and MP were partially supported by the Polish National Science Center under Ma\-estro Grant No. 2014/14/A/ST1/00453.
   KM was supported by NSF DMS-1521771, DMS-1622401, DMS-1839294,  DARPA contracts HR0011-16-2-0033 and FA8750-17-C-0054, and NIH grant R01 GM126555-01.
}

\begin{abstract}
The topological method for the reconstruction of dynamics from time series  \cite{MiMrReSz99}
is reshaped to improve its range of applicability, particularly in the presence of sparse data and strong expansion.
The improvement is based on a multivalued map representation of the data.
However, unlike the previous approach, it is not required that the representation has a continuous selector.
Instead of a selector, a recently developed new version of Conley index theory for multivalued maps \cite{BM2016,B2017}
is used in computations. The existence of a continuous, single-valued
generator of the relevant dynamics is guaranteed in the vicinity of the graph of the multivalued map constructed from data.
Some numerical examples based on time series derived from the iteration of H\'enon type maps are presented.
\end{abstract}
\maketitle
%-------------------------------------------------------------------

%!TEX root = ./ConSampDyn_2.tex

%-------------------------------------------------------------------------
%%%%%%%%%%%%%%%%%%%%%%%%%%%%%%%%%%%%%%%%%%%%%%%%%%%%%%%%%%%%%%%%%%%%%%%%%%%%%%%%%%%%%%%%%%%%
\section{Introduction}

Conceptual models for most physical systems are based on a continuum; values of the states of a system are assumed to be real numbers.
At the same time science is increasingly becoming data driven and thus based on finite information.
This suggests the need for tools that seamlessly and systematically provide information about continuous structures from finite data and accounts for the rapid rise in use of methods from topological data analysis (TDA).
However, not surprisingly, there are significant challenges associated with the sampling or generation of data versus the necessary coverage from which draw the appropriate conclusions.
In this paper we focus on this challenge in the context of nonlinear dynamics.

The fundamental work of  Niyogi, Smale, and Weinberger \cite{NSW2008} provides probabilistic guarantees that the correct homology groups have been computed, but is based on uniform sampling of the manifold.
For a nonlinear dynamical system one expects that the sampling is influenced by an underlying invariant measure that is rarely uniform with respect to the volume of the underlying phase space.
Furthermore, in practice one seldom knows the underlying subset of phase space on which the dynamics of interest occurs, e.g. the invariant set.
As a consequence one must expect that in applications we will need to collect considerably more data than a theoretical minimum would necessitate.

The predominant tool used by the TDA community to overcome the problem of lack of knowledge of the topological space of interest is persistent homology that provides homological information at all scales.
There are two challenges associated with this approach.
The first is that persistent homology computations on large data sets can be prohibitively expensive (there is extensive work being done to address this problem \cite{dey2013,oudot2015,harker2018}), and second that the development of a persistence theory of maps is in its early stages \cite{DJKKLM2019,EJM2015,BEJM2019}.
An alternative technique is to bin the data.
This is the approach we adopt in this paper.
In particular, we assume that the data points are measured via coordinates and thus the binning in phase space naturally takes the form of cubical sets.
The advantage is that we can a priori choose the bins so that the homological computations are feasible given time and memory constraints, and almost tautologically the binning process is a data reduction technique.

Identification of the space is only part of the challenge of understanding dynamics, we also need to capture the  behavior of the nonlinear map that generates the dynamics.
Though an oversimplification,  interesting dynamics is often driven by nonlinearities that exhibit significant expansion.
As is made explicit in \cite{FMN2014} the amount of data needed to expect a correct direct computation of the induced maps on homology is proportional to the magnitude of the Lipschitz constant of the map.
This will not be a surprise to anyone who has attempted to construct explicit simplicial maps for nonlinear functions.
The significance of the work reported in this paper is that we can obtain reliable information about the dynamics without directly identifying the map.

\begin{figure}[ht]
\begin{minipage}[t]{0.3\linewidth}
\centering
\includegraphics[width=\textwidth]{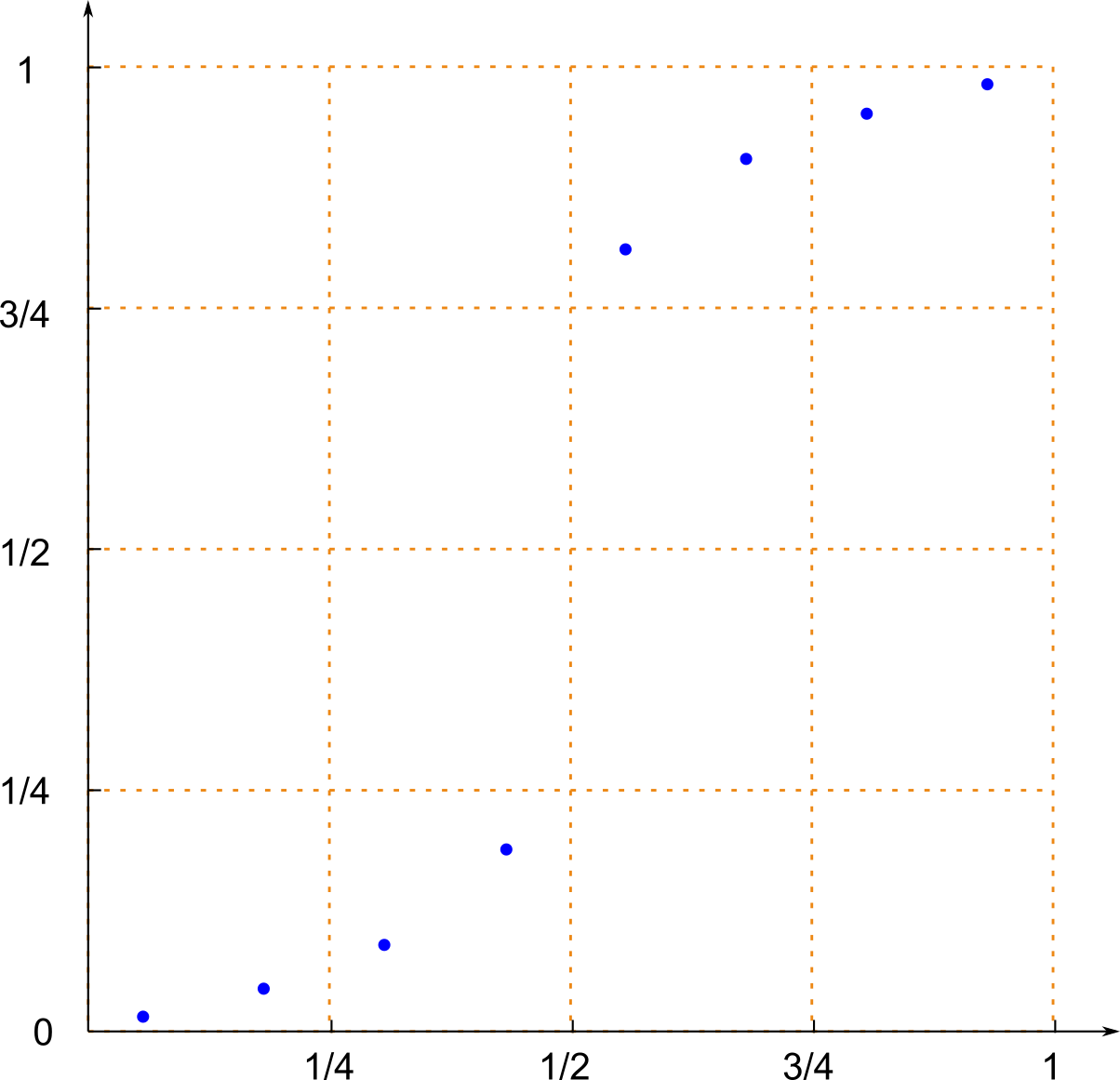}
\subcaption{The data marked by blue dots and the grid
indicated with orange dashed lines.}
\end{minipage}
\hspace{4mm}
\begin{minipage}[t]{0.3\linewidth}
\centering
\includegraphics[width=\textwidth]{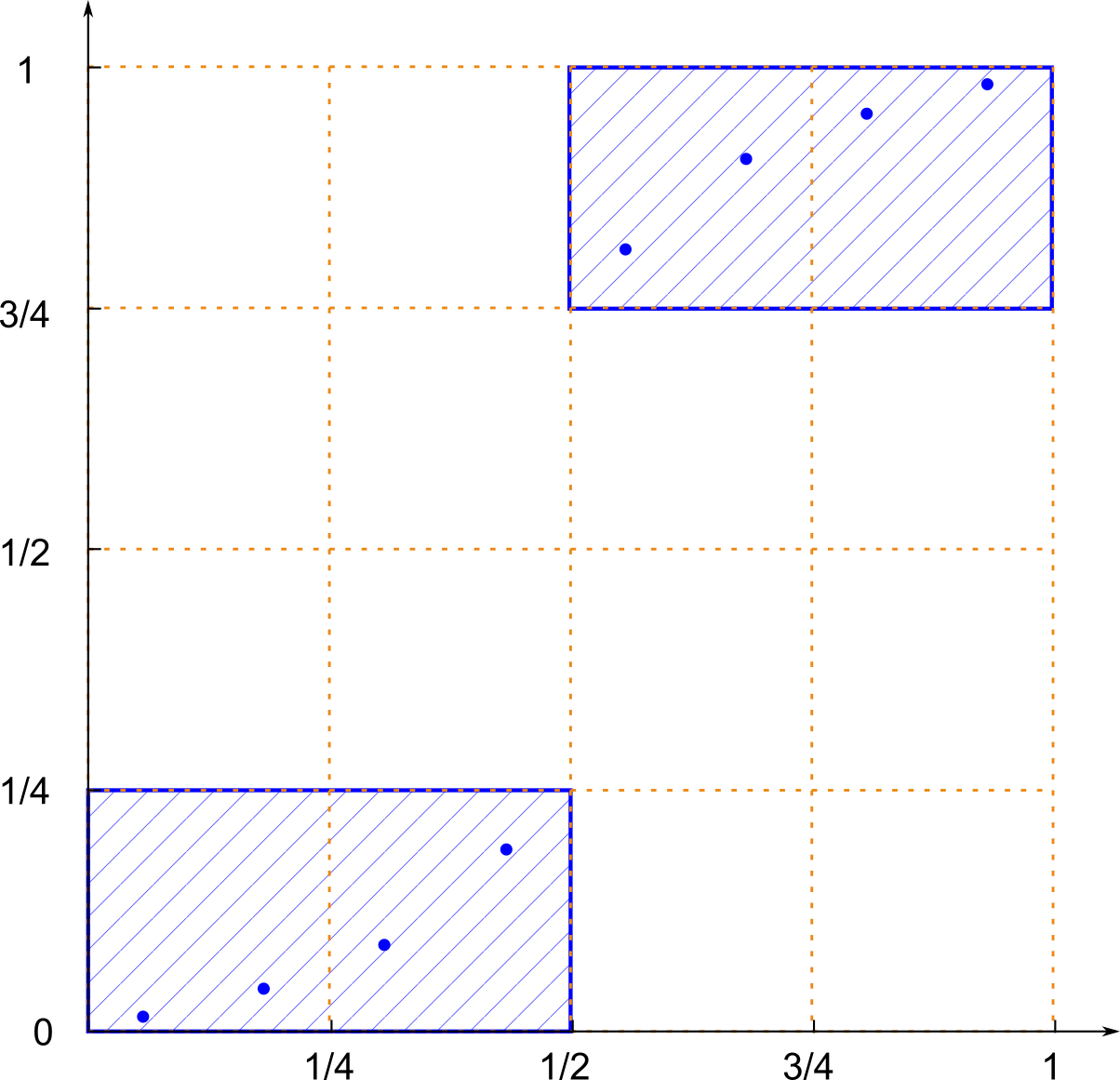}
\subcaption{The bins of data indicated with four squares shaded with blue.}
\end{minipage}
\hspace{4mm}
\begin{minipage}[t]{0.3\linewidth}
\centering
\includegraphics[width=\textwidth]{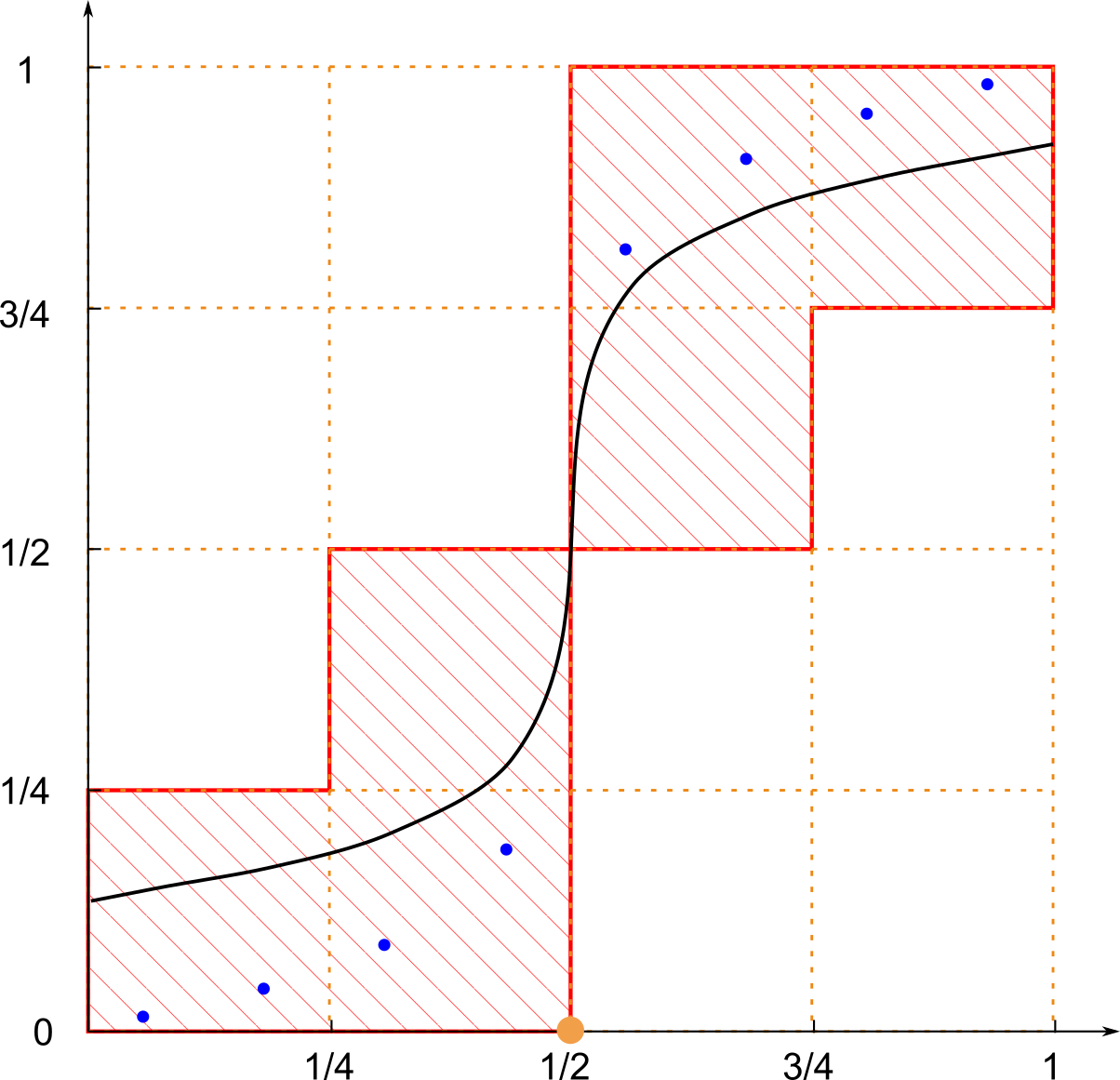}
\subcaption{The expansion of bins indicated with six squares shaded with red and the graph of a continuous selector in black.}
\end{minipage}\\
\begin{minipage}[t]{0.3\linewidth}
\centering
\includegraphics[width=\textwidth]{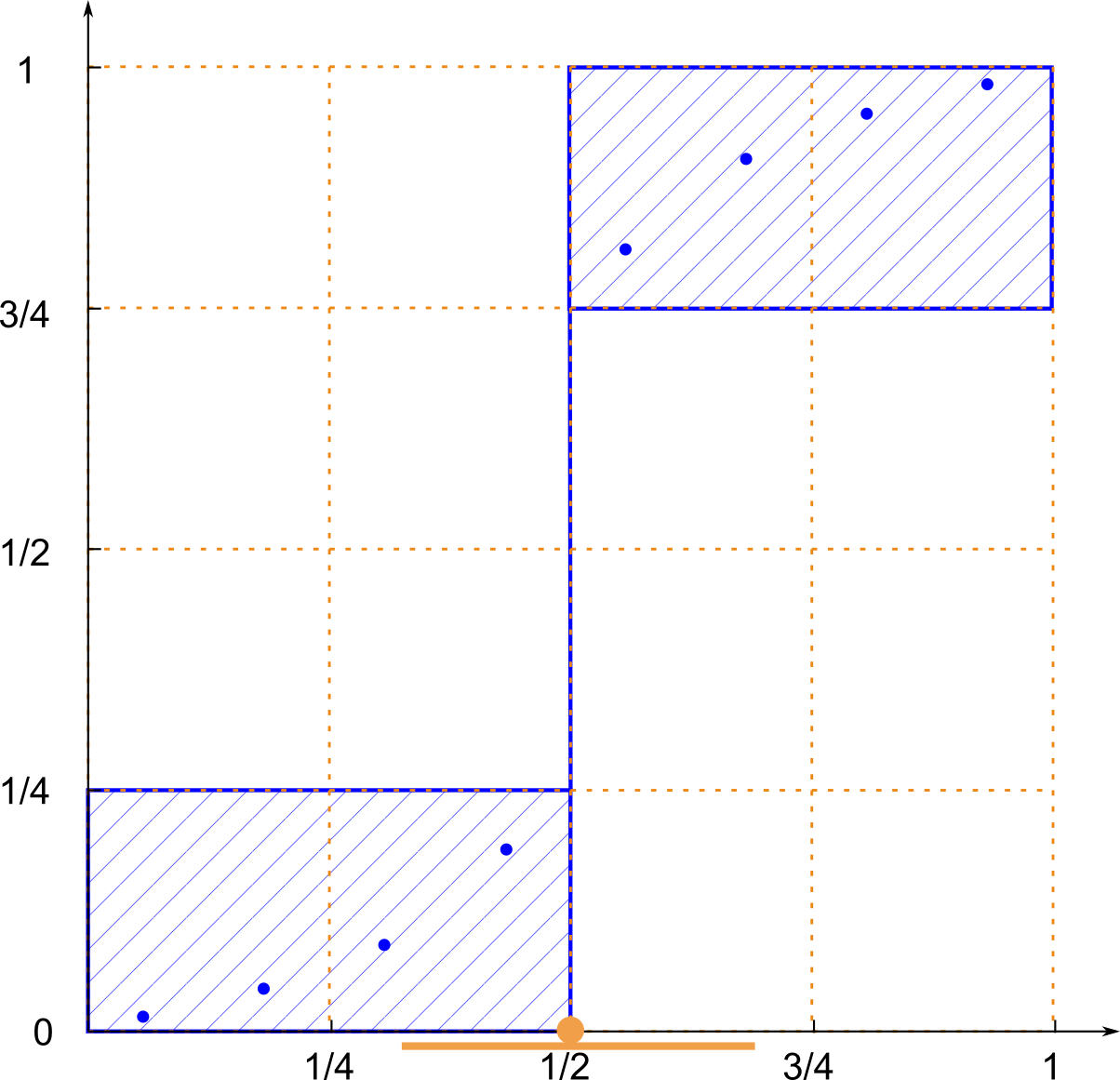}
\subcaption{The graph of an upper semicontinuous acyclic map $F:X\mto X$ in blue. Isolating neighborhood $N$ marked by orange line segmant.}
\end{minipage}
\hspace{5mm}
\begin{minipage}[t]{0.3\linewidth}
\centering
\includegraphics[width=\textwidth]{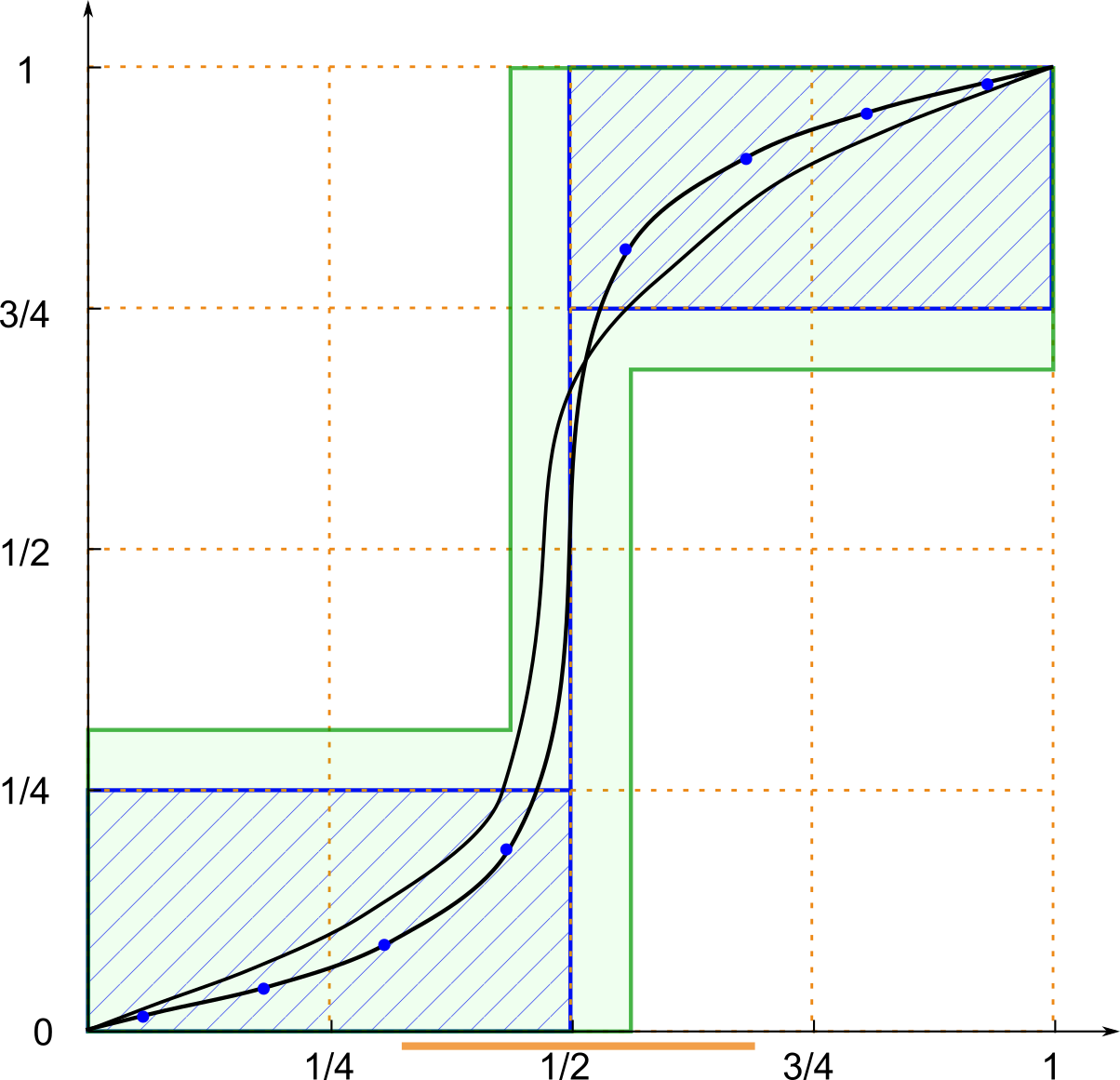}
\subcaption{The graph of map $F:X\mto X$ in blue, and its vicinity for continuous maps sharing with $F$ isolating neighborhood $N$ and the Conley index, in green.}
\end{minipage}
\caption{Construction of an upper semicontinuous acyclic multivalued map $F$ covering points representing the data.}
\label{fig:example}
\end{figure}

To explain the philosophy before becoming submerged in the technical details, consider a dynamical system on the unit interval and  assume that we have collected the data $\setof{(x,y)\in [0,1]\times [0,1] }$ as indicated in Figure~\ref{fig:example}(A).
We interpret this data as providing information about the graph of a continuous map $f\colon [0,1]\to [0,1]$ and
the question we ask is: \emph{can we extract information about the dynamics generated by $f$}?
The answer is yes.
In fact, under minimal hypotheses we can conclude that there are attractors that contain a fixed point within the intervals $[0,\frac{1}{4}]$ and
$[\frac{3}{4},1]$, and there exists an unstable invariant set, also containing a fixed point, in the interval $[\frac{3}{8}, \frac{5}{8}]$.
These results are obtained by building an upper semi-continuous acyclic multivalued map from the available data,
applying to it a recently developed new version of Conley index theory for multivalued maps \cite{BM2016,B2017}
in order to identify isolating neighborhoods and index pairs,
and then computing the associated Conley indices (definitions and details are provided in the following sections).
The last point requires that we be able to compute an induced map on homology.

An outline for the strategy used to perform these identifications is as follows.
As indicated above we bin the data.
Using intervals of length 1/4 to define the bins we obtain the blue shaded regions shown in Figure~\ref{fig:example}(B).
The blue regions are meant to provide a representation $F$ of the graph of the unknown function $f$.
Of course, as presented this is impossible; the domain of $F$ is connected but the blue regions are not.
One means of addressing this issue is to expand the representation so that the graph of a continuous functions can be included into the representation, i.e.\ the representation admits a \emph{continuous selector}.
Techniques of this type were successfully employed  in \cite{MiMrReSz99}.
However, they may easily fail. Applying the method of \cite{MiMrReSz99} to the representation in Figure~\ref{fig:example}(B)
leads to the representation in Figure~\ref{fig:example}(C). Actually, this is a minimal expansion which admits
a continuous selectors satisfying $f(\frac{1}{2}) = \frac{1}{2}$.
However, the resulting approximation of the dynamics is too crude: the combinatorial procedure for finding isolation neighborhoods presented in \cite{Szymczak-1997,Szymczak-1999} fails to produce an isolating neighborhood for the fixed point $x = \frac{1}{2}$.
On one hand, one can easily check that any other procedure must fail in this case, because the identity map is among selectors.
On the other hand, using an even larger expansion that produces an outer approximation
\cite{KMV2005} and using methods detailed in  \cite{D2003, DFT2008}
the desired isolating neighborhood and index pair can be recovered.
However, our experience is that applying this latter approach to complex time series data
even for two-dimensional examples often results in failure.

Lest the reader think that this is a contrived example, consider the function $f:[0,1]\to[0,1]$ given by
$f(x) = -nx^3+(1+n)x$
and observe that for $n\geq 1$ the points in Figure~\ref{fig:example} are consistent with data lying near the graph of $f$.
The dynamics generated by $f$ consists of stable fixed points at $0$ and  $1$, an unstable fixed point at $1/2$ and connecting orbits
from the unstable fixed point to the stable fixed points.
Furthermore, as $n$ increases, the minimal Lipschitz constant of $f$ given by $f'(1/2)=\frac{4+n}{4}$ increases
which results in dynamics becoming more pronounced.
However, from the perspective of experimental or numerically derived data, we expect the data points to cluster along the lines $y=0$ and $y=1$, and thus the observed discontinuity becomes more pronounced especially if one refines the binning.
We take this to be yet another suggestion that the direct approach of constructing a representation that admits a continuous selector is not the ideal technique.

As indicated above, we draw conclusions about the continuous dynamics from induced maps on homology via the Conley index.
This suggests that to obtain motivation for an alternative approach we consider the example from a purely  homological perspective.
Consider a function $f\colon [0,1]\to [0,1]$ and its graph $G_f:= \setof{(x,y)\in [0,1]^2\mid y=f(x)}$.
Let $\pi_1 \colon G_f \to  [0,1]$ and $\pi_2\colon G_f \to  [0,1]$ denote the projections from the graph to the domain and range of $f$, respectively.
Then $\pi_1$ is a homeomorphism, $\pi_{1*}$ is invertible and,
on the level of homology, $f_* = \pi_{2*}\circ \pi_{1*}^{-1}$.
Observe that if we replace $G_f$ by the blue shaded regions shown in Figure~\ref{fig:example}(B)
then $\pi_{1*}$ is not invertible, but we still can deduce
the correct map  induced by $F$ on homology. This is because the pre-image  $\pi_{1*}^{-1}$ takes on two values,
but these values are mapped to the same value under $\pi_{2*}$.
For a more complete discussion on this perspective see \cite{HKMP2016}.
What should be clear is that to apply this in general we require a condition that forces $\pi_{2*}$ to collapse appropriate generators in
the homology of the representation $H_*(F)$.

With this in mind consider the blue region shown in Figure~\ref{fig:example}(D).
In this case the fiber's of $\pi_1$ are acyclic, thus  $\pi_{1*}$ is invertible, and the question of how $\pi_{2*}$ acts on generators is resolved.
Because we are interested in extracting dynamics, rather than considering the blue region to be a fiber bundle over the phase space, we view it as the graph of an upper semi-continuous acyclic multivalued map $F: [0,1]\mto [0,1]$ and we use $F$ to extract isolating neighborhoods, index pairs, and ultimately the Conley index.

We note that in this simple one-dimensional example, the choice of the blue line in Figure~\ref{fig:example}(D) is obvious.
In higher dimensions there are a variety of means of attempting to resolve the issue of controlling how $\pi_{2*}$ acts on generators from the pre-image of $\pi_{1*}$ and the identification of optimal methods remains an open question.
In this paper we seek minimal rectangular regions.

To be more specific we assume that our data consists of a finite set of points $A\subset \RR^d$ and our understanding of the dynamics is to be derived from the map $g\colon A\to \RR^d$.
We also assume that we have chosen a scale $\delta >0$ for the binning and that the bins take the form \[
[n_1\delta, (n_1+1)\delta] \times [n_2\delta, (n_2+1)\delta] \times \cdots \times [n_d\delta, (n_d+1)\delta]
\]
where $n_i\in\ZZ$.
More generally, we work with \emph{$\delta$-cuboids}, sets of the form
\begin{equation}\label{df:delta_cuboid}
[n_1\delta,m_1\delta]\times[n_2\delta,m_2\delta]\times\cdots\times[n_d\delta,m_d\delta],
\end{equation}
where  $(n_1,n_2,\ldots,n_d),(m_1,m_2,\ldots,m_d)\in\ZZ^d$.
An {\em elementary cube} is a cuboid where $m_i-n_i \in\{0,1\}$ for $i=1,2,\ldots d$.
We denote the set of all $\delta$-cuboids in $\RR^d$ by $\cC^d_\delta$ and the set of all $\delta$-cubes in $\RR^d$ by $\cK^d_\delta$.

For a bounded subset $X\subset\RR^d$ we introduce the following notation.
\[
  K_\delta(X):=\bigcup\setof{Q\in\cK^d_\delta\mid X\cap Q\neq\emptyset},
\]
  and
 \[
  \llcorner X \urcorner_\delta :=\bigcup\setof{Q\in\cK^d_\delta\mid \conv{(X)}\cap Q\neq\emptyset},
\]
where $\conv{(X)}$ denotes the convex hull of $X$.

Returning to the map  $g\colon A\to \RR^d$  its  {\em sunflower enclosure} is the multivalued map $F^s _{g,\delta}\colon  K_\delta(A)\mto\RR^d$ defined by
\[
   F^s _{g,\delta}(x):=\llcorner g(K_\delta(x)\cap A)\urcorner_\delta \subset \RR^d .
\]
We leave it to the reader to check that given $\setof{(x,g(x))\in [0,1]\times [0,1] }$ as shown in Figure~\ref{fig:example}(A),  the graph of $F^s _{g,\delta}$ is as shown in  Figure~\ref{fig:example}(D).

Sunflower enclosures satisfy a variety of nice properties.
Recall (cf. \cite{KMM}) that $F\colon X\mto \RR^d$ is {\em cubical} if
\begin{itemize}
\item[(a)] $X\subset \RR^n$ is a \emph{cubical} set, i.e. it can be written as a finite union of elementary cubes,
\item[(b)] for any $x\in X$ the set $F(x)$ is cubical,
\item[(c)] for any elementary cube $Q = [n_1\delta,m_1\delta]\times\cdots\times[n_d\delta,m_d\delta]$ in $X$, $F_{|\mathring{Q}}$ is constant, where
$\mathring{Q} :=  (n_1\delta,m_1\delta)\times\cdots\times(n_d\delta,m_d\delta)$ and
$(n_i\delta,m_i\delta) = \{n_i\}$ if $n_i=m_i$.
\end{itemize}
The following proposition follows from \cite[Proposition 14.5]{G76}.
\begin{prop}
\label{prop:sunflower-enclosure}
 A sunflower enclosure is an upper semicontinuous cubical map.
 \end{prop}
When the values of the sunflower enclosure are contractible, then using algorithms developed in \cite{Szymczak-1997} and the formula from \cite[Theorem 4.4]{B2017} one can identify cubical isolating blocks,  cubical weak index pairs and an index map associated with  $F^s _{g,\delta}$ (see \cite{Pr19} for more details).
In particular, a Conley index can be computed.

From the perspective of identifying dynamics  the aforementioned computation should be viewed as purely formal, e.g.\ in and of itself it does not guarantee that there is a continuous map that generates dynamics that is compatible with the associated Conley indices.
The majority of this paper is dedicated to guaranteeing that the formal computation does in fact lead to the existence of a large, but explicit, family of nonlinearities that are capable of producing the observed dynamics.
To state our goals more precisely we introduce the following notation.
Let $F\colon X\mto X$.
For simplicity of notation we identify $F$ with its graph $\setof{(x,y)\in X\times X\mid y\in F(x)}$.
Using the max-norm on the product space $X\times X$, let $B(F,\varepsilon)\subset X\times X$ denote the open set of points within $\varepsilon$ of the graph of $F$ (see  Figure~\ref{fig:example}(E)).
Following \cite{GGK} (cf. e.g. \cite{G76}) we say that a continuous single valued map $f:X\to X$ is a {\em continuous $\varepsilon$-approximation (on the graph) of $F:X\mto X$} if $f\subset B(F,\varepsilon)$.

We denote the set of continuous $\varepsilon$-approximations of $F$ by $a_\varepsilon(F)$.

Our claim is that Conley index information computed for $F\colon X\mto X$, an acyclic upper semicontinous cubical maps, is valid for the dynamics generated by any continuous function $f\in a_\varepsilon(F)$ for all $\varepsilon\in (0,\varepsilon_0)$ sufficiently small.
As the  results described below indicate, our approach provides explicit lower bounds on $\varepsilon_0$.

We have, up to this point in the introduction, be rather circumspect about how the Conley index provides information about nonlinear dynamics.
One of the more powerful results is that it can be used to construct semi-conjugacies to known dynamics.
To be more precise, given two continuous maps $f\colon X\to X$ and $\sigma\colon Y\to Y$, $f$ is \emph{semi-conjugate} to $\sigma$ if there exists a continuous surjective map $\rho\colon X\to Y$ such that
\[
\begin{tikzcd}
X\ar{d}[swap]{\rho}\ar{r}{f}&X\ar{d}{\rho}\\[3ex]
Y\ar{r}[swap]{\sigma}&Y
\end{tikzcd}
\]
commutes.
Semi-conjugacies are of interest if the dynamics of $\sigma$ is understood, as this implies that the dynamics of $f$ must be at least as complicated, i.e.\ one can deduce structure about the dynamics of $f$ from that of $\sigma$.

In the context of the Conley theory, one begins with an index pair $P = (P_1,P_2)$ (see Section~\ref{sec:preliminary} for precise definitions).  The homological Conley index is derived from a map  $f_{P*} \colon H_*(P_1/P_2, [P_2]) \to H_*(P_1/P_2, [P_2])$ that itself is derived from the action of $f$ on the pointed quotient space $(P_1/P_2, [P_2])$.
Let $N = \cl(P_1\setminus P_2)$.  The meta form of the desired theorem is that given the homological Conley index,  information about the index pair, and an explicit dynamical system $\sigma\colon Y\to Y$, then there exists a semi-conjugacy
\[
\begin{tikzcd}
\Inv(N,f)\ar{d}[swap]{\rho}\ar{r}{f}& \Inv(N,f) \ar{d}{\rho}\\[3ex]
Y\ar{r}[swap]{\sigma}&Y
\end{tikzcd}
\]
where $\Inv(N,f)$ denotes the maximal invariant set in $N$ under $f$.

The potential of the proposed theory in applications is demonstrated in \cite{BMMP2018app},
in particular in examples based on the time series studied in \cite{MiMrReSz99}.
In this paper we will prove the following two results.

\begin{figure}
\begin{center}
  \includegraphics[width=0.8\textwidth]{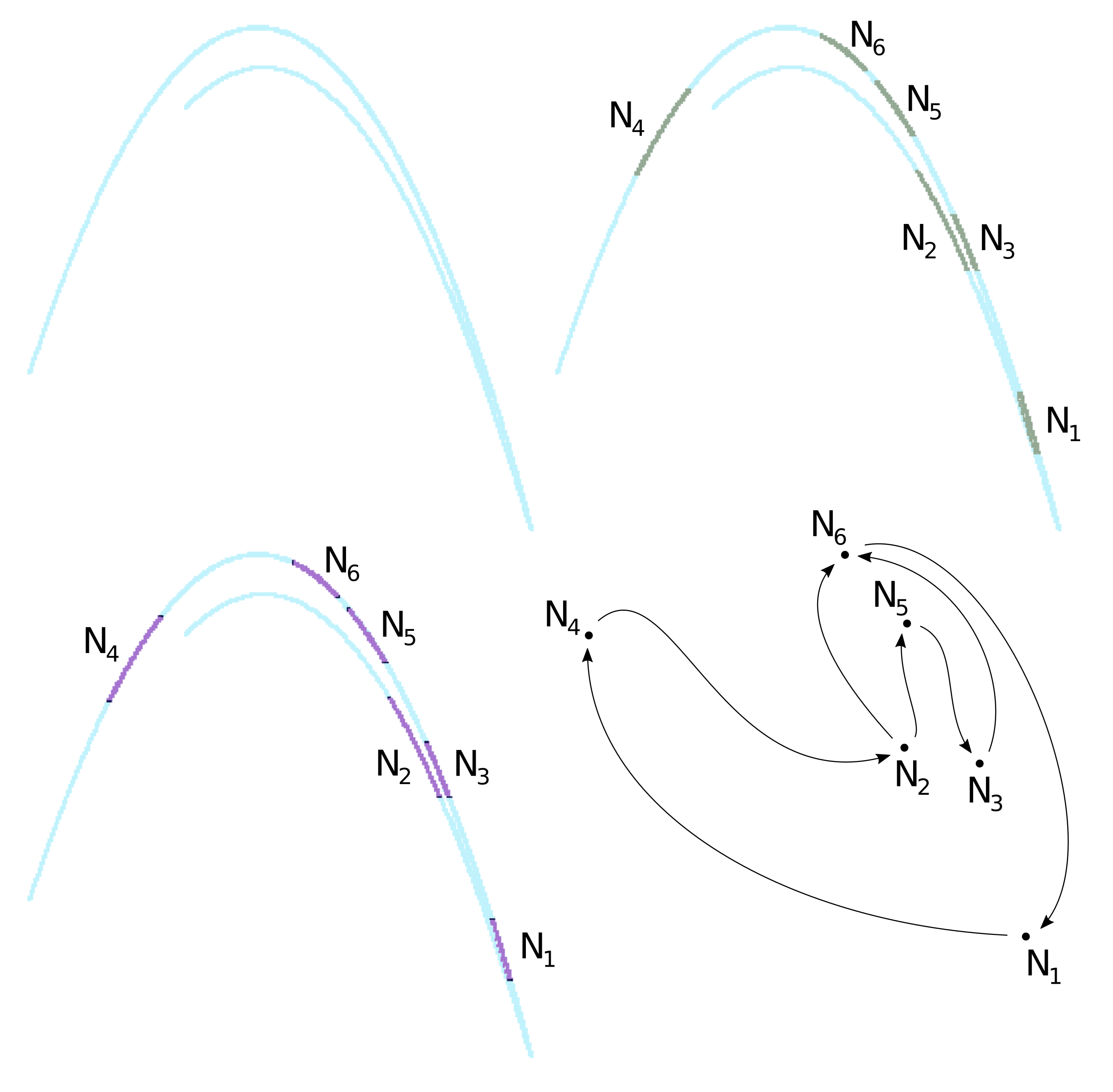}
\end{center}
  \caption{Domain of sunflower enclosure for $g_{\bar{x}}$, an isolating neighbourhood (in dark sea green), its weak index pair (in blue violet) and the graph of transitions between components of an isolating neighborhood.}
  \label{fig:henon2d}
\end{figure}

\begin{thm}
\label{thm:henon2d}
Consider the time series $\bar{x}=(x_i)_{i=100}^{30000}$ generated by iterating the H\'enon map
\[
     H:\RR^2\ni (x,y) \mapsto (1 - a x^2 + b y, x) \in\RR^2
\]
 with the parameter values $a = 1.65$, $b = 0.1$, and initial condition $(x_0,y_0)=(0,0)$.
 Set
\[
A_{\bar{x}} := \setof{ (x_i, x_{i+1}) \mid i=100,\ldots, 29,999 }
\]
and let $g_{\bar{x}} \colon A_{\bar{x}}\to \RR^2$ be given by $g_{\bar{x}}(x_i,x_{i+1}) = (x_{i+1},x_{i+2})$.

Choose a binning of $\RR^2$ based on $\delta:=0.008127$ and let $F := F^s _{g_{\bar{x}},\delta}\colon  K_\delta(A_{\bar{x}})\mto\RR^2$ be the sunflower enclosure of $g_{\bar{x}}$, i.e.
\[
   F^s _{g_{\bar{x}},\delta}(x):=\llcorner g(K_\delta(x)\cap A_{\bar{x}})\urcorner_\delta \subset \RR^2.
\]
Let $\varepsilon = \delta/2$.

Then, $a_\varepsilon(F)\neq\emptyset$.
Furthermore, there exists a compact set $N\subset\RR^2$ (see Figure~\ref{fig:henon2d}) such that for any $f\in a_\varepsilon(F)$
\begin{itemize}
   \item[(i)] $N$ is an isolating neighborhood of $f$,
   \item[(ii)] there exists a semiconjugacy $\theta_f: \Inv(N,f)\to \Sigma_A$ onto the subshift of finite type on six symbols with the transition
   matrix
\[
    A=\left(
\begin{array}{cccccc}
0 & 0 & 0 & 0 & 0 & 1 \\
0 & 0 & 0 & 1 & 0 & 0 \\
0 & 0 & 0 & 0 & 1 & 0 \\
1 & 0 & 0 & 0 & 0 & 0 \\
0 & 1 & 0 & 0 & 0 & 0 \\
0 & 1 & 1 & 0 & 0 & 0 \\
\end{array}
\right)
\]
such that for every periodic $a\in\Sigma_A$ there exists a periodic point of $f$ in $\theta_f^{-1}(a)$.
\end{itemize}
In particular, $f$ has positive topological entropy on $\Inv(N,f)$.
\end{thm}

\begin{figure}
\begin{center}
  \includegraphics[width=0.8\textwidth]{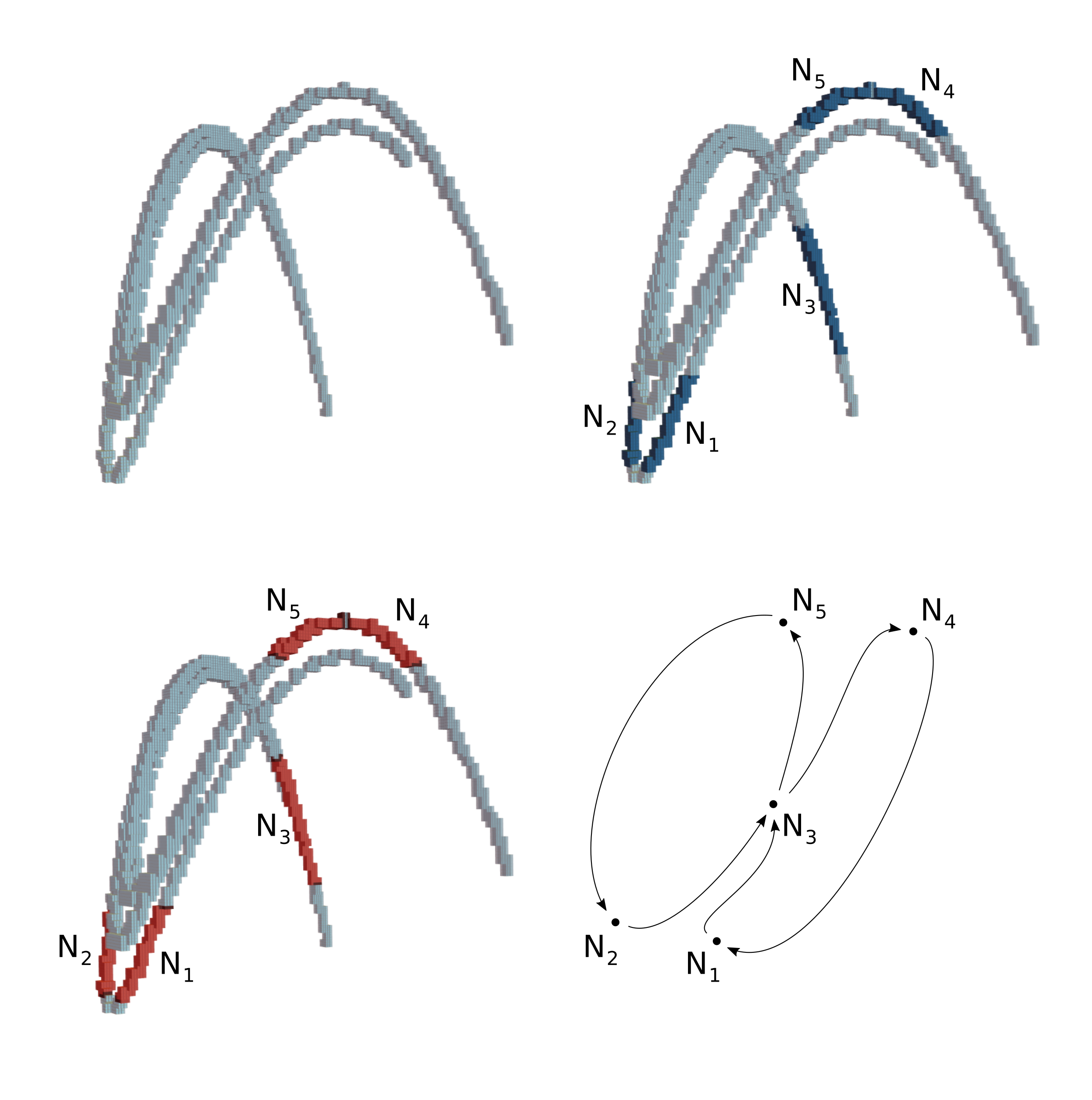}
\end{center}
  \caption{Domain of sunflower enclosure for $g_{\bar{x}}$, an isolating neighbourhood (in dark cyan), its weak index pair (in orange) and the graph of transitions between components of an isolating neighborhood.}
  \label{fig:henon3d}
\end{figure}

\begin{thm}
\label{thm:henon3d}
Consider the time series $\bar{x}=(x_i)_{i=100}^{30000}$ generated by iterating the delayed H\'enon map
\[
     H:\RR^3\ni (x,y,z) \mapsto (1 - a x^2 + b z, x, y) \in\RR^3
\]
with the parameter values $a = 1.65$, $b = 0.1$, and initial point $(x_0,y_0,z_0) = (0,0,0)$.
Set
\[
A_{\bar{x}} := \setof{ (x_i, x_{i+1}, x_{i+2}) \mid i=100,\ldots, 29,998 }
\]
and let $g_{\bar{x}} \colon A_{\bar{x}}\to \RR^3$ be given by $g_{\bar{x}}(x_i,x_{i+1},x_{i+2}) = (x_{i+1},x_{i+2},x_{i+3})$.

Choose a binning of $\RR^3$ based on $\delta:= 0.035258$ and let $F := F^s _{g_{\bar{x}},\delta}\colon  K_\delta(A_{\bar{x}})\mto\RR^3$ be the sunflower enclosure of $g_{\bar{x}}$.

Let $\varepsilon = \delta/2$.

Then, $a_\varepsilon(F)\neq\emptyset$.
Furthermore, there exists a compact set $N\subset\RR^3$ (see Figure~\ref{fig:henon3d}) such that for any $f\in a_\varepsilon(F)$
\begin{itemize}
   \item[(i)] $N$ is an isolating neighborhood of $f$, and
   \item[(ii)] there exists a semiconjugacy $\theta_f: \Inv(N,f)\to \Sigma_A$ onto the subshift of finite type on five symbols with the transition matrix
\[
    A=\left(
\begin{array}{ccccc}
 0 & 0 & 0 & 1 & 0  \\
 0 & 0 & 0 & 0 & 1  \\
 1 & 1 & 0 & 0 & 0  \\
 0 & 0 & 1 & 0 & 0  \\
 0 & 0 & 1 & 0 & 0  \\
\end{array}
\right)
\]
such that for every periodic $a\in\Sigma_A$ there exists a periodic point of $f$ in $\theta_f^{-1}(a)$.
\end{itemize}
In particular, $f$ has positive topological entropy on $\Inv(N,f)$.
\end{thm}

A natural question arises how sensitive these results are to the choice of $\delta$, the length of the time series, or the choice of initial condition.
The fundamental feature of the Conley index is that it does not change under a small perturbation of the generator of the dynamical systems.
Thus, the question reduces to the understanding of the stability of the multivalued map representation of the data.
It is natural to expect that by increasing the length of the time series or changing the initial condition the semiconjugacy should be preserved
as long as the same isolating neighborhood is used. Experiments we run confirm this expectation.
A more delicate is the question how the choice of $\delta$ affects the results.
On one hand, if $\delta$ is very small, then the domain of the multivalued representation becomes a collection of isolated cubes.
Therefore, it cannot properly approximate the phase space which is a continuum.
On the ohter hand, if  $\delta$ is too large, the  multivalued representation gives a very coarse description of dynamics.
Therefore, one cannot expect that it will give an interesting description of dynamics.
Thus, the optimum is somewhere in the middle. Experiments we run show that small changes to $\delta$ preserve the results and
moderate changes lead to a different matrix $A$ but still let us claim the existence of an invariant set with positive entropy.
An interesting problem is to get the understanding of changes in the results under varying $\delta$ in the spirit of persistent homology.
This is left for future investigations.

We now provide an outline for the paper.
Section~\ref{sec:preliminary} provides basic definitions related to the Conley index.
Section~\ref{sec:uppersemi} presents results about isolating neighborhoods in the context of upper semi-continuous multivalued maps.
Section~\ref{sec:approx} makes use of the results of Section~\ref{sec:uppersemi} to provide conditions under which continuous functions in a neighborhoods of the graph of a upper semi-continuous mulitvalued map $F$ with convex compact images inherit isolating neighborhoods and their associated Conley index from $F$.
Results of this form are essential.
The isolating neighborhood and Conley index computations in Theorems~\ref{thm:henon2d} and \ref{thm:henon3d} are done using the sunflower enclosure $F$, but the results of interest concern the dynamics generated by continuous in $a_{\varepsilon}(F)$.

The conclusion of Theorems~\ref{thm:henon2d} and \ref{thm:henon3d}  involve the existence of a semi-conjugacy.
As indicated above this is done via the Conley index.
Because we work with upper semi-continuous multivalued maps that need not admit a continuous selector, we need to work with \emph{weak index pairs}.
The classical result of Szymczak \cite{Szymczak-1996,Szymczak-1997} that proves the existence of a semi-conjugacy onto symbolic dynamics is based on a stronger definition of an index pair and therefore cannot be applied directly.
Section~\ref{sec:semiconjugacy} presents theorems that are an extension of Szymczak's results.
Sections~\ref{sec:indexmap} - \ref{sec:lefschetz} provide the necessary background to prove the results of Section~\ref{sec:semiconjugacy}.

The fact that Theorems~\ref{thm:henon2d} and \ref{thm:henon3d}  contains explicit bounds on the class of maps, e.g.\ $a_\varepsilon(F)$ with $\varepsilon = \delta/2$ is important for the development of models.
Section~\ref{sec:hve} provides explicit information about the preservation of topological and dynamical properties for continuous functions near $F$.

Finally, the proofs of Theorems~\ref{thm:henon2d} and \ref{thm:henon3d} are presented in Section~\ref{sec:henon}.

%%%%%%%%%%%%%%%%%%%%%%%%%%%%%%%%%%%%%%%%%%%%%%%%%%%%%%%%%%%%%%%%%%%%%%%%%%%%%%%%%%%%%%%%%%%%
\section{Preliminaries}
\label{sec:preliminary}

Throughout this paper by an interval in the set of integers $\ZZ$ we mean the intersection of a closed interval in $\RR$ with $\ZZ$. For $n\geq 1$ let $I_n:=\{1,2,\ldots,n\}$ and for $p\geq 2$ let $\ZZ_p:=\{0,1,\ldots,p-1\}$ denote the additive topological group with addition modulo $p$ and discrete topology.

Given a topological space $X$ and a subset $A\subset X$, by $\Int_X A$, $\cl_X A$ will denote the {\em interior} of $A$ in $X$ and the {\em closure} of $A$ in $X$ respectively. We omit the symbol of space if the space is clear from the context. 

Let $X$, $Y$ be topological spaces. By $F\colon X\mto Y$ we denote a multivalued map, that is a map $F\colon X\ni x\mapsto F(x)\in \cP(Y)$, where $\cP(Y)$ is the power set of $Y$. A multivalued map $F$ is {\em upper semicontinuous} if for any closed $B\subset Y$ its large counter image under $F$, that is the set $F^{-1}(B):=\{x\in X\ |\ F(x)\cap B\neq\emptyset\}$, is closed. 

Throughout the paper we identify $F$ with its {\em graph}, that is the set $\{(x,y)\in X\times Y\ |\ y\in F(x)\}$. In the sequel, we are interested in multivalued self-maps, that is multivalued maps of the form $F\colon X\mto X$.

Let $I$ be an interval in $\ZZ$ with $0\in I$. A single valued mapping $\sigma\colon I\to X$ is a \emph{solution for $F$ through $x\in X$} if $\sigma(n+1)\in F(\sigma(n))$ for all $n, n+1\in I$ and $\sigma(0)=x$ (cf. \cite[Definition 2.3]{KM95}).

Given a subset $N\subset X$, the set 
$$
\Inv(N,F):=\{x\in N\ |\ \exists\sigma\colon\ZZ\to N \text{ a solution for } F \text{ through } x\}
$$ 
is called the {\em invariant part} of $N$. A compact subset $N\subset X$ is an {\em isolating neighborhood} for $F$ if $\Inv(N,F)\subset\Int N$. A compact subset $N\subset X$ is called an {\em isolating block} with respect to $F$ if
$$
N\cap F(N)\cap F^{-1}(N)\subset \Int N.
$$
Note that any isolating block is an isolating neighborhood. A compact set $S\subset X$ is said to be {\em invariant} with respect to $F$ if $S = \Inv(S,F)$. It is called an {\em isolated invariant set} if it admits an isolating neighborhood $N$ for $F$ such that $S = \Inv(N,F)$ (cf. \cite[Definition 4.1, Definition 4.3]{BM2016}). 

By $F$-{\em boundary} of a given set $A\subset X$ we mean $\bd_{F}A := \cl A \cap \cl(F(A)\setminus A)$. Let $N\subset X$ be an isolating neighborhood for $F$. 

\begin{defn}[cf. {\cite[Definition 4.7]{BM2016}}]
\label{defn:wip}
A pair $P=(P_1,P_2)$ of compact sets $P_2\subset P_1\subset N$ is called a {\em weak index pair} in $N$ if
\begin{itemize}
\item[(a)] $F(P_i)\cap N\subset P_i$ for $i\in\{1,2\}$,
\item[(b)] $\bd_F P_1\subset P_2$,
\item[(c)] $\Inv(N,F)\subset \Int(P_1\setminus P_2)$,
\item[(d)] $P_1\setminus P_2\subset \Int N$.
\end{itemize}
\end{defn}

Given a weak index pair $P$ in an isolating neighborhood $N\subset X$ for $F$ we set
$$
T_N(P):=(T_{N,1}(P),T_{N,2}(P)):=(P_1\cup (X\setminus\Int N),P_2\cup (X\setminus\Int N)).
$$
Recall (cf. e.g. \cite{BM2016,Mr06}) that $F_P$, the restriction of $F$ to the domain $P$, is a multivalued map of pairs, $F_P\colon P\mto T_N(P)$; the inclusion $i_P:P\to T_N(P)$ induces an isomorphism in the Alexander--Spanier cohomology; and the {\em index map} $I_{F_P}$ is defined as an endomorphism of $H^*(P)$ given by
$$
I_{F_P}=F_P ^*\circ (i_P ^*)^{-1}.
$$
The pair $(H^*(P),I_{F_P})$ is a graded module equipped with an endomorphism. Applying the Leray functor $L$ (cf. \cite{M90,BM2016}) to $(H^*(P),I_{F_P})$ we obtain a graded module with its endomorphism which we call the {\em Leray reduction of the Alexander--Spanier cohomology of a weak index pair} $P$.

\begin{defn}[cf. {\cite[Definition 6.3]{BM2016}}]
\label{defn:conley_index}
The graded module $L(H^{*}(P),I_{F_P})$, that is the Leray reduction of the Alexander--Spanier cohomology of a weak index pair $P$ is called the {\em cohomological Conley index of} $\Inv(N,F)$ and denoted by $C(\Inv(N,F),F)$. 
\end{defn}

%----------------------------
%%%%%%%%%%%%%%%%%%%%%%%%%%%%%%%%%%%%%%%%%%%%%%%%%%%%%%%%%%%%%%%%%%%%%%%%%%%%%%%%%%%%%%%%%%%%
\section{Dynamics of upper semicontinuous maps}
\label{sec:uppersemi}
Let $(X,d)$ be a metric space. By $B_r(x)$ we denote the open ball with the center in $x\in X$ and radius $r>0$. Closed balls will be denoted by $\bar{B_r}(x)$. For a given $A\subset X$,  $B_r(A)$ will stand for an {\em open} {\em $r$-hull of $A$}, that is,
$$
B_r(A):=\bigcup\{B_r(a)\,|\,a\in A\}.
$$

Let $F\colon X\mto X$ be an upper semicontinuous map. One can easily verify that (multivalued) selections of $F$ share with $F$ its isolating neighborhood and a weak index pair. We express this observation here for further references.
\begin{prop}
\label{prop:wip_for_mv_selector}
Assume $N$ is an isolating neighborhood for an upper semicontinuous $F\colon X\mto X$, $P$ is a weak index pair for $F$ in $N$ and $G\colon X\mto X$ is an upper semicontinuous map such that $G\subset F$. Then $N$ is an isolating neighborhood for $G$, and  $P$ is a weak index pair for $G$ in $N$.
\end{prop}
%%%
The aim of this section is to show that, to a certain extent, the reverse implications holds true. To be precise, we have the following theorem.
\begin{thm}\label{thm:NFepsilon}
Let  $N$ be an isolating neighborhood with respect to an upper semicontinuous map $F:X\mto X$. There exists an $\varepsilon >0$ such that $N$ is an isolating neighborhood with respect to an arbitrary upper semicontinuous map $G:X\mto X$ with $G\subset B(F,\varepsilon)$.
\end{thm}
We postpone its proof to the end of this section.

\begin{lem}\label{lm:l1}
Let $A\subset X$ be a compact set and let $\{x_n\}\subset X$ be a sequence convergent to $x\in X$. If $x_n\in B(A,\frac{1}{n})$ for $n\in\NN$ then $x\in A$.
\end{lem}
\begin{proof}
Suppose the contrary and consider an $r>0$ such that $B(x,r)\cap A=\emptyset$. Observe that, for large enough $n\in\NN$, we have $d(x_n,x)\leq \frac{r}{2}$. Moreover, there exists a sequence $\{u_n\}\subset A$ with $d(u_n,x_n)\leq \frac{1}{n}$ for $n\in\NN$. However, $d(x_n,u_n)\geq d(u_n,x)-d(x_n,x)\geq r-\frac{r}{2}=\frac{r}{2}$, a contradiction.
\qed\end{proof}
%-----------------------
\begin{lem}\label{lm:l2}
Let $F:X\mto X$ be upper semicontinuous and let $N\subset X$ be compact. A solution $\tau:\ZZ\to N$ for $F$ through $x\in N$ exists provided for any $n\in \NN$ there exists a solution $\sigma:[-n,n]\to N$ through $x$.
\end{lem}
\proof
Let $\sigma ^n:[-n,n]\to N$ be a solution with respect to $F$ through $x$. By induction we construct a sequence of solutions $\tau ^n:[-n,n]\to N$ for $F$ through $x$ such that
\begin{itemize}
\item[(p1)] there exists a strictly increasing sequence $\{m_p\}\subset\NN$ such that $\tau ^n (k)=\lim_{p\to\infty}\sigma ^{m_p}(k)$ for any $k\in [-n,n]$,
\item[(p2)] $\tau ^{n-1}\subset \tau^n$ for $n\geq 1$.
\end{itemize}
Define $\tau ^0:[0]\to N$ by putting $\tau ^0(0):=x$. Clearly (p1) and (p2) hold. Suppose $\tau ^n$ has been constructed so that (p1) and (p2) hold. Denote $\bar{\sigma} ^p:=\sigma ^{m_p}$ and take into account a subsequence $\bar{p}$ such that the sequences $\bar{\sigma} ^{\bar{p}}(n+1)$ and $\bar{\sigma} ^{\bar{p}}(-n-1)$ converge to $v,w\in N$, respectively. We define $\tau ^{n+1}:[-n-1,n+1]\to N$ by
$$
\tau ^{n+1}(k):=\left\{\begin{array}{rl}
\tau ^n (k)&\for |k|\leq n\\
v&\for k=n+1\\
w&\for k=-n-1.
\end{array}\right.
$$
It is straightforward to see that conditions (p1) and (p2) hold, and $\tau^{n+1}(0)=x$. It remains to be verified that $\tau ^{n+1}$ is a solution for $F$. Since $\bar{\sigma} ^{p}$ is a solution for $F$, we have
\begin{equation}\label{eq:rsigma}
\bar{\sigma} ^{\bar{p}}(k+1)\in F(\bar{\sigma} ^{\bar{p}}(k))\for k\in\ZZ.
\end{equation}
For any $k\in [-n-1,n+1]$ the sequence $\bar{\sigma} ^{\bar{p}}(k)$ converges to $\tau^{n+1}(k)$. Because the graph of $F$ is closed (cf. \cite[Proposition 14.4]{G76}), passing to the limit in (\ref{eq:rsigma}) we have $\tau^{n+1}(k+1)\in F(\tau^{n+1}(k))$.
\qed\\
%--------------

{\bf Proof of Theorem \ref{thm:NFepsilon}}. 
For contradiction suppose that for any $m\in\NN$ there exists an upper semicontinuous $G_m:X\mto X$ with $G_m\subset B(F,\frac{1}{m})$ and such that $\Inv (N, G_m)\cap\bd N\neq \emptyset$. Let $x_m\in \Inv (N, G_m)\cap\bd N$. Passing to a subsequence, if necessary, we may assume that $x_m$ converges to an $x\in\bd N$. Let $\sigma _m:\ZZ\to N$ be a solution for $G_m$ through $x_m$. Fix an integer $n\in\NN$, choose a subsequence $m_p$ such that for any $k\in [-n,n]$ the sequence $\sigma _{m_p}(k)$ is convergent, and define $\tau ^n:[-n,n]\to N$ by putting $\tau ^n (k):=\lim_{p\to \infty}\sigma _{m_p}(k)$ for $k\in [-n,n]$. We have $(\sigma _{m_p}(k),\sigma _{m_p}(k+1))\in G_{m_p}\subset B(F,\frac{1}{m_p})$. Using Lemma \ref{lm:l1} we infer that $(\tau ^n(k),\tau ^n(k+1))\in F$, which means that $\tau ^n:[-n,n]\to N$ is a solution for $F$ through $x$. This, along with Lemma \ref{lm:l2}, yields the existence of a solution $\tau:\ZZ\to N$ for $F$ through $x$. However, $x\in\bd N$, a contradiction.
\qed

%%%%%%%%%%%%%%%%%%%%%%%%%%%%%%%%%%%%%%%%%%%%%%%%
\section{$\varepsilon$-Approximations}\label{sec:approx}
In the following we consider the Cartesian product of normed spaces as the normed space with the max-norm.

Following \cite{GGK} (cf. e.g. \cite{G76}) we say that a continuous single valued map $f:X\to X$ is a {\em continuous $\varepsilon$-approximation (on the graph) of $F:X\mto X$} if $f\subset B_\varepsilon(F)$. We denote the set of continuous $\varepsilon$-approximations of $F$ by $a_\varepsilon(F)$.

%-----------------------
\begin{thm}\label{thm:homotop_star}
Let $Y$ be a normed space and let $X\subset Y$ be compact. Assume that $F:X\mto X$ is an upper semicontinuous map with convex and compact values, and $N$ is an isolating neighborhood with respect to $F$. Then:
\begin{itemize}
\item[(i)] there exists an $\varepsilon_0 > 0$ such that, for any $0<\varepsilon\leq\varepsilon_0$ there is  a continuous $\varepsilon$-approximation $f:X\to X$ of $F$ such that $N$ is an isolating neighborhood with respect to $f$, and $C(\Inv(N,F),F)=C(\Inv(N,f),f)$;
\item[(ii)] if $X$ is an ANR then there exists a $\delta>0$ such that for any continuous $\delta$-approximation $g:X\to X$ of $F$ we have $C(\Inv(N,F),F)=C(\Inv(N,g),g)$.
\end{itemize}
\end{thm}
\begin{proof}
Take an $\varepsilon_0>0$ as in Theorem  \ref{thm:NFepsilon} and $0<\varepsilon\leq\varepsilon_0$. By \cite[Theorem 1]{C69}
there exists a continuous $\varepsilon$-approximation $f:X\to X$ of $F$. We shall prove that $f$ satisfies the assertions (i) and (ii).

To this end, for $\lambda\in [0,1]$, we define $F_\lambda:X\mto X$ by
$$
F_\lambda (x):=\lambda f(x)+(1-\lambda)F(x)\for x\in X.
$$
It follows from  the upper semicontinuity of $F$ and the continuity of $f$ that $F_\lambda$ is upper semicontinuous and it is straightforward to observe that $F_\lambda $ has convex and compact values.

According to the construction of the $\varepsilon$-approximation $f$ of $F$ in \cite{C69}, for arbitrarily fixed $x\in X$ there exists an $x'\in B_\varepsilon(x)$ such that $f(x)\in B_\varepsilon(F(x'))$ and $F(x)\subset B_\varepsilon(F(x'))$.
Therefore, for any $\lambda\in [0,1]$, we have $F_\lambda(x)\subset B_\varepsilon(F(x'))$, as $B_\varepsilon(F(x'))$ is convex. Consequently, $F_\lambda\subset B_\varepsilon(F)$.
Theorem \ref{thm:NFepsilon} shows that $N$ is an isolating neighborhood with respect to $F_\lambda$ for every $\lambda\in[0,1]$. Therefore, by the continuation property of the Conley index (cf. \cite[Theorem 6.1]{B2017}), we have $C(\Inv(N,F),F)=C(\Inv(N,f),f)$.

Let an $\varepsilon >0$ be as above. By \cite[Theorem 23.9]{G76} there is a $\delta \in (0,\varepsilon]$ such that for any $f,g:X\to X$, the $\delta$-approximations of $F$, there exists a homotopy $h:X\times [0,1]\to X$ joining $f$ and $g$, such that $h(\cdot, t)$ is an $\varepsilon$-approximation of $F$, for all $t\in [0,1]$.
Fix such a $\delta>0$ and consider $f:X\to X$, a $\delta$-approximation of $F$ defined as in \cite{C69}. Let $g:X\to X$ be an arbitrary $\delta$-approximation of $F$. Since $\delta\leq\varepsilon$, by \cite[Theorem 5.13]{GGK} (cf. e.g. \cite[Theorem 23.9]{G76}) and Theorem  \ref{thm:NFepsilon}, $\Inv(N,f)$ and $\Inv(N,g)$ are related by continuation; hence $C(\Inv(N,g),g)=C(\Inv(N,f),f)$. This, along with property (i), completes the proof.
\qed\end{proof}
%----------------------------
%%%%%%%%%%%%%%%%%%%%%%%%%%%%%%%%%%%%%%%%%%%%%%%%%%%%%%%%%%%%%%%%%%%%%%%%%%%%%%%%%%%%%%%%%%%%
\section{$\varepsilon$-Approximations of cubical maps}
\label{sec:hve}

In this section we assume that $X\subset\RR^d$ is a closed subset and $F\colon X\mto X$ is a multivalued cubical map (cf. e.g. \cite{KMM}), and $\varrho$ stands for the max metric in $\RR^d$. For $x\in\RR^d$ by $\sigma_x$ we denote the unique elementary cube such that $x\in\mathring{\sigma}_x$. For $\varepsilon>0$ define maps $F_{\varepsilon}, F^{\varepsilon}\colon X \mto X$ by
\begin{equation}
\label{eq:F_sub_epsilon}
F_{\varepsilon}(x):= F(\Be(x))
\end{equation}
and
\begin{equation}
\label{eq:F_sup_epsilon}
F^{\varepsilon}(x):=\Be(F(x)).
\end{equation}
We refer to maps $F_\varepsilon$ and $F^\varepsilon$ as a {\em horizontal}  and  a {\em vertical enclosure} of $F$, respectively. 

We begin with some auxiliary lemmas.

\begin{lem}
\label{lem:distance_in_intersections}
Assume $A_1, A_2\subset X$ are cubical, $\varepsilon\in(0,\frac{1}{2})$ and $y\in\Be(A_1)\cap\Be(A_2)$. Then, there exists a $y'\in A_1\cap A_2$ such that $\varrho(y,y')\leq 2\varepsilon$.
\end{lem}
\proof
For $i=1,2$ let $y_i\in A_i$ be such that $\varrho(y_i,y)<\varepsilon$. Then $\sigma_{y_1}\cap\sigma_{y_2}\neq\emptyset$ and $\sigma_{y_i}\subset A_i$. Let $y'\in\sigma_{y_1}\cap\sigma_{y_2}$. Then $y'\in A_1\cap A_2$ and $\varrho(y,y')\leq\varrho(y,y_1)+\varrho(y_1,y')\leq2\varepsilon$.
\qed

\begin{lem}
\label{lem:inclusion_induces_isomorphism}
Assume $P\subset M\subset\RR^d$ are cubical and $0<\varepsilon < \frac{1}{2}$. Then the inclusion $\mu\colon P\cup M\to \Be(P)\cup M$ induces isomorphism in cohomology. 
\end{lem}
%------------
\proof
Consider the multivalued map $G\colon \Be(P)\cup M\to M$ given by
$$
G(x):=\{y\in M\ |\ \varrho(x,y) = \varrho(x,M)\}.
$$
This map has compact values and is upper semicontinuous (see \cite[Lemma 1]{Mr91}). Since $G(x) = \{x\}$ for $x\in M$, we see that $G\circ\mu = \id_{P\cup M}$. We will show that $\mu\circ G$ is homotopic to $\id_{\Be(P)\cup M}$. One easily verifies that $G(x)=\bar{B}_{\varrho(x,M)}(x)\cap M$. Let
$$
\cQ:=\{Q\in\cK\ |\ Q\subset M, Q\cap\bar{B}_{\varrho(x,M)}(x)\neq\emptyset\}
$$
and $\cQ':=\{Q\cap\Be(x)\ |\ Q\in\cQ\}$. Then $G(x) = \bigcup_{Q'\in\cQ'}Q'$. Each $Q'\in\cQ'$ is a rectangle as an intersection of rectangles. Hence it is convex. We claim that $\bigcap_{Q'\in\cQ'}Q'\neq\emptyset$. For this end it suffices to show that $Q_1'\cap Q_2'\neq\emptyset$ for any  $Q_1', Q_2'\in\cQ'$. Since $Q_i'\neq\emptyset$, $i\in\{1,2\}$, take $x_i\in Q_i'$. Then $x_1, x_2\in \Be(x)$, which means that $\varrho(x_1,x_2)\leq 2\varepsilon$. Since $\varepsilon<\frac{1}{2}$, this implies that $Q_1'\cap Q_2'\neq\emptyset$. Thus, $F(x)$ is star-shaped, hence acyclic. 

For $\lambda\in[0,1]$ let
$$
G_\lambda(x):=\{(1-\lambda)x+\lambda y\ |\ y\in G(x)\}
$$
and $D(x):=\bigcup_{\lambda\in[0,1]}G_\lambda(x)$. Note that if $x\in M$ then $G_\lambda(x)=D(x)=\{x\}$. Also if $x\in\Be(P)\cup M$ then $D(x)\subset\Be(P)\cup M$. Therefore, 
$$
[0,1]\times (\Be(P)\cup M)\ni (\lambda,x)\longmapsto G_\lambda(x)\subset \Be(P)\cup M
$$
is the requested homotopy between $\mu\circ G$ and $\id_{\Be(P)\cup M}$.
\qed

As a consequence of the previous lemma we have the following lemma.
\begin{lem}
\label{lem:homot_equiv}
Assume $A\subset\RR^d$ is a cubical set and $0<\varepsilon<\frac{1}{2}$. Then $A$ and $ \Be(A)$ are homotopy equivalent.  \qed
\end{lem}

Now we enumerate a few properties of the enclosures.
\begin{lem}
\label{lem:F_sub_epsilon}
The map $F_\varepsilon$ has the following properties:
\begin{itemize}
\item[(i)] If $A\subset X$, then $F_{\varepsilon}^{-1}(A) = \Be(F^{-1}(A))$.
\item[(ii)] If $F$ is upper semicontinuous, then so is $F_\varepsilon$.
\item[(iii)] If $\varepsilon<\frac{1}{2}$ and $F$ is upper semicontinuous, then for any $x\in X$ there is $y\in X$ with $F(y)=F_{\varepsilon}(x)$.
\item[(iv)] If $\varepsilon<\frac{1}{2}$ and $F$ is upper semicontinuous and convex-valued, then so is $F_\varepsilon$.
\item[(v)] If $\varepsilon<\frac{1}{2}$ and $F$ is upper semicontinuous and has contractible values, then so does $F_\varepsilon$.
\item[(vi)] If $\varepsilon<\frac{1}{2}$ and $F$ is an upper semicontinuous map with convex values then $F_\varepsilon$ admits a continuous selection.
\item[(vii)] If $A\subset X$ is a cubical set and $F$ is upper semicontinuous, then $F_\varepsilon(A) = F(A)$ whenever $\varepsilon<\frac{1}{2}$.
\end{itemize}
\end{lem}
\proof
In order to prove inclusion $\Be(F^{-1}(A)) \subset F_{\varepsilon}^{-1}(A)$ in (i) take ant $x\in\Be(F^{-1}(A))$ and an $x'\in F^{-1}(A)$ such that $x\in\Be(x')$. Then $F(x')\cap A \neq\emptyset$. Take a $y\in F(x')\cap A$. Then $y\in F_\varepsilon(x)$ and $F_\varepsilon(x)\cap A\neq\emptyset$ which proves that $x\in F_\varepsilon^{-1}(A)$. 

In reverse direction, take an $x\in F_\varepsilon^{-1}(A)$, a $y\in F_\varepsilon(x)\cap A$ and an $x'\in\Be(x)$ such that $y\in F(x')$. It means that $F(x')\cap A\neq\emptyset$ and $x'\in F^{-1}(A)$. Therefore $x\in\Be(F^{-1}(A))$.

By (i), the large counterimage under $F_\varepsilon$ of any closed set in its range is closed. Hence, $F_\varepsilon$ is upper semicontinuous, and we have (ii) .

In order to show (iii), fix $x\in X$ and consider the set
$$
\cA (x):=\left\{\mathring{Q}\;|\;\mathring{Q}\mbox{ is a cell}, \mathring{Q}\cap\bar{B_\varepsilon}(x)\neq\emptyset\right\}.
$$
Note that $\cA (x)\neq\emptyset$. Since $\varepsilon<\frac{1}{2}$, for any $\mathring{Q},\mathring{Q'}\in \cA (x)$ we have $Q\cap Q'\neq\emptyset$. Then $P:=\bigcap _{\mathring{Q}\in \cA (x)}Q\neq\emptyset$ and $P$ is an elementary cube. Moreover, $P$ is a face of every cube $Q$ with $\mathring{Q}\in \cA (x)$. Then, for any $y\in\mathring{P}$ we have $F(y)=F(\bar{B_\varepsilon}(x))=F_\varepsilon (x)$, as $F$ is upper semicontinuous.

Properties (iv) and (v) follow from (iii).

We shall prove (vi). To this end consider $\tilde{F}_\varepsilon:X\mto X$ given by
$$
\tilde{F}_\varepsilon(x):=F(B_\varepsilon(x))\for x\in X.
$$
It is easy to see, by the same reasoning as for $F_\varepsilon$, that $\tilde{F}_\varepsilon$ has nonempty convex and compact values.
Moreover, the large counterimage under $\tilde{F}_\varepsilon$ of any single point in its range is open; hence $\tilde{F}_\varepsilon$ is lower semicontinuous. Consequently, by the Michael's selection theorem (cf. e.g. \cite{Mi56}), there exists a continuous map $f:X\to X$, a selection of $\tilde{F}_\varepsilon$. Clearly $\tilde{F}_\varepsilon\subset F_\varepsilon$, hence $f$ is a continuous selection of $F_\varepsilon$, as desired.

In order to prove (vii) take a $y\in F_\varepsilon(A)$. Then there exist an $x\in A$ such that $y\in F_\varepsilon(x)$ and an $x'\in\bar{B}(x,\varepsilon)$ such that $y\in F(x')$. Then $\sigma_x\cap\sigma_{x'}\neq \emptyset$, because $\varrho(x,x')\leq 2\varepsilon< 1$. It follows that we can take an $x''\in\sigma_x\cap\sigma_{x'}$. By the strong upper semicontinuity of $F$ we set $F(x')\subset F(x'')$. Hence, $y\in F(x'')\subset F(\sigma_x) \subset F(A)$. The inclusion in the reverse direction is obvious. 
\qed

\begin{lem}
\label{lem:F_sup_epsilon}
The map $F^\varepsilon$ has the following properties:
\begin{itemize}
\item[(i)] If $A\subset X$, then $F^{\varepsilon}(A) = \Be(F(A))$.
\item[(ii)] If $F$ is usc, then so is $F^\varepsilon$.
\item[(iii)] If $F$ is convex-valued, then so is $F^\varepsilon$.
\item[(iv)] If $\varepsilon<\frac{1}{2}$  and a cubical map $F$ has contractible values, then so does $F^\varepsilon$.
\item[(v)] If $A\subset X$ is a cubical set and $F$ is cubical, then $(F^\varepsilon)^{-1}(A) = F^{-1}(A)$ for any $0\leq\varepsilon<1$.
\end{itemize}
\end{lem}
\proof
To prove (i) observe that
$$
\begin{array}{rcl}
F^\varepsilon(A) & = & \bigcup_{x\in A} F^\varepsilon(x) =  \bigcup_{x\in A}\Be(F(x))  = \\
& = & \bigcup_{x\in A}\bigcup_{y\in F(x)} \Be(y) =  \bigcup_{y\in F(A)} \Be(y)  = \Be(F(A)). \\
\end{array}
$$

Properties (ii) and (iii) are obvious. 

Property (iv) is a consequence of Lemma \ref{lem:homot_equiv}.

In order to show inclusion $F^{-1}(A)\subset (F^\varepsilon)^{-1}(A)$ in (v) take $x\in F^{-1}(A)$. It means that $F(x)\cap A \neq \emptyset$ and $F^\varepsilon(x)\cap A\neq\emptyset$. Hence $x\in (F^\varepsilon)^{-1}(A)$. 

To prove the opposite inclusion take an $x\in(F^\varepsilon)^{-1}(A)$. Since $F^\varepsilon(x)\cap A\neq \emptyset$, there exist a $y\in  F^\varepsilon (x)\cap A$ and a $y'\in F(x)$ such that $y\in\Be(y')$. We have $\sigma_y\cap\sigma_{y'}\neq\emptyset$, because $\varrho(y,y')\leq\varepsilon<1$. Take $y''\in\sigma_y\cap\sigma_{y'}$. Then $y''\in\sigma_{y'} = \cl\mathring{\sigma}_{y'}\subset F(x)$, because $F$ is cubical. Notice that $y\in\mathring{\sigma}_y\cap A$. Thus, $y''\in\sigma_y\subset A$, because $A$ is a cubical set. It follows that $F(x)\cap A\neq\emptyset$ and $x\in F^{-1}(A)$ which completes the proof.   
\qed

Horizontal enclosures preserve isolating neighborhoods and weak index pairs. More precisely, we have the following propositions. 

%--------------
\begin{prop}
\label{thm:preserving_iso_ne}
Assume $F\colon X\mto X$ is a cubical, upper semicontinuous multivalued map and $N$ is a cubical isolating neighborhood for $F$. Then, for any $\varepsilon<1$, we have $\Inv(N,F_\varepsilon)\subset\Be(\Inv(N,F))$. As a consequence, $N$ is an isolating neighborhood for $F_\varepsilon$.
\end{prop}
\proof
Take an arbitrary $x_0\in\Inv(N,F_\varepsilon)$ and consider $x\colon\ZZ\to N$, a solution for $F_\varepsilon$ in $N$ through $x_0$. Let $n\in\ZZ$ be fixed. We have $x_{n+1}\in F_\varepsilon(x_n)$. There exists an $x'_n\in N$ such that $F_\varepsilon(x_n) = F(x'_n)$ and $\varrho(x'_n,x_n)\leq\varepsilon<1$. Therefore we can take an $x''_n\in\sigma_{x_n}\cap\sigma_{x'_n}$ such that $\varrho(x''_n,x_n)\leq\varepsilon$ and $F(x'_n)\subset F(x''_n)$, as $F$ is upper semicontinuous. We have $x_{n+1}\in F(x'_n)$ and $\sigma_{x_{n+1}}\subset F(x'_n)$, because $F$ is cubical. Moreover, $x''_{n+1}\in F(x'_n)\subset F(x''_n)$ and $\varrho(x''_n,x_n)\leq\varrho(x'_n,x_n)\leq\varepsilon$. Since $n\in\ZZ$ was arbitrarily fixed, we have constructed $x''\colon\ZZ\to N$, a solution for $F$ in $N$, with $\varrho(x''_n,x_n)\leq\varepsilon$. In particular, $x _0\in \Be(\Inv(N,F))$, showing that $\Inv(N,F_\varepsilon)\subset \Be(\Inv(N,F))$.

Since $\Inv(N,F)\subset\Int N$ and $\varepsilon<1$, the latter inclusion yields $\Inv(N,F_\varepsilon)\subset\Be(\Inv(N,F))\subset \Int N$. This completes the proof.
\qed 
%------------

\begin{prop}
\label{thm:preserving_wip}
Assume $F\colon X\mto X$ is a cubical, upper semicontinuous multivalued map, $N$ is a cubical isolating neighborhood for $F$, $P$ is a cubical weak index pair in $N$, and $\varepsilon<\frac{1}{2}$.
Then $P$ is a weak index pair for $F_\varepsilon$. 
\end{prop}
\proof 
Properties (a) and (b) are straightforward consequences of Lemma \ref{lem:F_sub_epsilon}(vii). 

By Theorem \ref{thm:preserving_iso_ne} we have $\Inv(N,F_\varepsilon)\subset \Be(\Inv(N,F))\subset \Int(P_1\setminus P_2)$, and property (c) follows.

Property (d) is obvious.
\qed

%%%%%%%%%%%%%%%%%%%%%
Note that, in general, $F_\varepsilon$ is not a cubical map. However, it inherits from $F$ the following property.
\begin{lem}\label{lem:F_eps_prop}
If $F\colon X\mto X$ is a cubical map and $\varepsilon<\frac{1}{2}$, then
\begin{equation}\label{eq:F_eps_prop}
F_\varepsilon(y)\subset F_\varepsilon(x)\mbox{ whenever }\sigma_x\subset\sigma_y.
\end{equation}
\end{lem}
\proof
Since $\varepsilon<\frac{1}{2}$ and $\sigma_x\subset\sigma_y$, for an  arbitrary elementary cube $\sigma$, condition $\sigma\cap \bar{B}_\varepsilon (y)\neq\emptyset$ implies $\sigma\cap \bar{B}_\varepsilon (x)\neq\emptyset$. Therefore, taking into account that $F$ is cubical, we have 
$$
\begin{array}{rcl}
F_\varepsilon (y)&=&F(\bar{B}_\varepsilon (y))\\
&=&\bigcup_{z\in \bar{B}_\varepsilon (y)}F(z)\\
&=&\bigcup_{\sigma\cap \bar{B}_\varepsilon (y)\neq\emptyset}F(\mathring{\sigma})\\
&\subset& \bigcup_{\sigma\cap \bar{B}_\varepsilon (x)\neq\emptyset}  F(\mathring{\sigma})\\
&=&\bigcup_{z\in \bar{B}_\varepsilon (x)}F(z)\\
&=&F(\bar{B}_\varepsilon (x))\\
&=&F_\varepsilon (x).
\end{array}$$ 
This completes the proof.
\qed
%--------------------
\begin{prop}
\label{prop:preserving_iso_ne2}
Assume $F\colon X\mto X$ is a cubical, upper semicontinuous multivalued map and $N$ is a cubical isolating neighborhood for $F$. Then $N$ is an isolating neighborhood for $(F_\varepsilon)^\varepsilon$ for any $\varepsilon<\frac{1}{2}$.
\end{prop}
\proof 
For contradiction, suppose that $x_0\in\bd N$ and $x\colon\ZZ\to N$ is a solution for $(F_\varepsilon)^\varepsilon$ in $N$ through $x_0$. 

Let $n\in\ZZ$ be fixed. We have $x_{n+1}\in (F_\varepsilon)^\varepsilon (x_n)$. There exists an $x'_{n+1}\in F_\varepsilon (x_n)$ with $\varrho(x'_n,x_n)\leq\varepsilon<\frac{1}{2}$. Then $\sigma_{x_{n+1}}\cap\sigma_{x'_{n+1}}\neq\emptyset$, and we can take $x''_{n+1}\in\sigma_{x_{n+1}}\cap\sigma_{x'_{n+1}}$. Since $N$ is cubical and $x_{n+1}\in N$, we infer that $\sigma_{x_{n+1}}\subset N$. Hence, $x''_{n+1}\in\sigma_{x_{n+1}}\cap\sigma_{x'_{n+1}}\subset N$. Similarly, $x''_{n+1}\in\sigma_{x_{n+1}}\cap\sigma_{x'_{n+1}}\subset F_\varepsilon (x_n)$, as $x'_{n+1}\in F_\varepsilon (x_n)$ and $F_\varepsilon (x_n)$ is a cubical set. By Lemma \ref{lem:F_eps_prop} we have $F_\varepsilon (x_n)\subset F_\varepsilon (x''_n)$. This, along with $x''_{n+1}\in F_\varepsilon (x_n)$ yields $x''_{n+1}\in F_\varepsilon (x''_n)$. Since $n\in\ZZ$ was arbitrarily fixed, we have defined $x'':\ZZ\to N$, a solution with respect to $F_\varepsilon$ in $N$.

Note that $x''_0\in \bd N$, because $x_0\in \bd N$, $x''_0\in \sigma_{x_0} $ and $\bd N$ is a cubical set. This contradicts Proposition \ref{thm:preserving_iso_ne}, and completes the proof. 
\qed

%%%%%%%%%%%%%%%%%%%%%%%%%%%%%%%%%%%%%%%%%%%%%%
The following theorem is a counterpart of Theorem \ref{thm:homotop_star} for cubical maps.
%-----------------------
\begin{thm}\label{thm:homotop_star_2}
Assume that $F:X\mto X$ is an upper semicontinuous cubical map with contractible values and $\varepsilon <\frac{1}{2}$. Then $a_\varepsilon(F)\neq\emptyset$. Moreover, if $N$ is a cubical isolating neighborhood with respect to $F$, then $N$ is an isolating neighborhood with respect to arbitrary $f\in a_\varepsilon(F)$, and $C(\Inv(N,f),f)=C(\Inv(N,F),F)$.
\end{thm}
\proof The existence of an $\varepsilon$-approximation $f:X\to X$ of $F$ follows from \cite[Theorem 5.12]{GGK} (cf. e.g. \cite[Theorem 23.8]{G76}).

Since $F$ has contractible values then, by Lemma \ref{lem:F_sub_epsilon}(v) and Lemma \ref{lem:F_sup_epsilon}(iv), so does $(F_\varepsilon)^\varepsilon$. Moreover, by Proposition \ref{prop:preserving_iso_ne2}, $N$ is an isolating neighborhood with respect to $(F_\varepsilon)^\varepsilon$. Therefore we have well defined Conley index $C(\Inv(N,(F_\varepsilon)^\varepsilon),(F_\varepsilon)^\varepsilon)$. Since, in addition $F\subset (F_\varepsilon)^\varepsilon$, we infer that $C(\Inv(N,F),F)=C(\Inv(N,(F_\varepsilon)^\varepsilon),(F_\varepsilon)^\varepsilon)$. Note that if $f:X\to X$ is an $\varepsilon$-approximation of $F$, then we have $f\subset (F_\varepsilon)^\varepsilon$, and the identity $C(\Inv(N,f),f)=C(\Inv(N,(F_\varepsilon)^\varepsilon),(F_\varepsilon)^\varepsilon)$ follows. This completes the proof.
\qed
%-------

A statement analogous to Proposition \ref{thm:preserving_wip} for map $F^\varepsilon$ is not true, however an approximate version holds.
\begin{thm}
\label{thm:wip_for_vert_enclosure}
Assume that $F\colon X\mto X$ is an upper semicontinuous map with cubical values satisfying (\ref{eq:F_eps_prop}), $P$ is a cubical weak index pair with respect to $F$ in a cubical isolating neighborhood $N$, and $\varepsilon<\frac{1}{2}$. Then $\Be(P)$ is a weak index pair for $F^\varepsilon$ in $\Be(N)$.
\end{thm}
\proof
For the proof of property (a) fix an $i\in\{1,2\}$ and take an $x\in\Be(P_i)$ and a $y\in F^\varepsilon(x)\cap\Be(N)$. Then, there exists an $x'\in P_i$ such that $\varrho(x,x')<\varepsilon$ and by Lemma \ref{lem:distance_in_intersections} there exists a $y'\in F(x)\cap N$ such that $\varrho(y,y')<2\varepsilon<1$. Then $\sigma_x\cap\sigma_{x'}\neq\emptyset$ and we can take an $x''\in\sigma_x\cap\sigma_{x'}$ such that $\varrho(x,x'')<\varepsilon$. By (\ref{eq:F_eps_prop}) we have $F(x'')\supset F(x)$ and since $P$ is cubical, we have $x''\in\tau\subset P_i$. Similarly, $\varrho(y,y')<2\varepsilon<1$ implies that there exists a $y''\in \sigma_y\cap\sigma_{y'}$ such that $\varrho(y,y'')<\varepsilon$. Since $F(x)$ and $N$ are cubical and $y'\in\mathring{\sigma}_{y'}\cap F(x)\cap N$, we get $\sigma_{y'}\subset F(x)\cap N$. Therefore $y''\in F(x)\cap N\subset F(x'')\cap N$. By property (a) of $P$, we have $y''\in P_i$. Hence, $y\in\Be(P_i)$. 

In order to prove property (b) assume the contrary. Let $x\in\bd_{F^\varepsilon}\Be(P_1)\setminus \Be(P_2)$. It means that $x\in\Be(P_1)$, $x\in\cl(F^\varepsilon(\Be(P_1))\setminus\Be(P_1))$ and $x\notin \Be(P_2)$. Take an $x'\in P_1$ such that $\varrho(x,x')<\varepsilon$. Then $x'\notin P_2$, that is $x'\in P_1\setminus P_2$. Consider a sequence $(x_n)_{n\in\NN}$ such that $x_n\in F^\varepsilon(\Be(P_1))\setminus\Be(P_1)$ and $x_n\to x$. It follows that for every $n\in\NN$ we have $x_n\in F^\varepsilon(u_n)$ for some $u_n\in\Be(P_1)$. Take a $u_n'\in P_1$ such that $\varrho(u_n,u_n')<\varepsilon$ and a $z_n\in F(u_n)$ such that $\varrho(x_n,z_n)<\varepsilon$. We have $z_n\notin P_1$, because otherwise $x_n\in\Be(P_1)$. By $\varrho(u_n,u'_n)<\varepsilon$, we can take $u''_n\in\sigma_{u_n}\cap\sigma_{u'_n}$. Since $P$ is cubical, we have $u''_n\in P_1$. Since $F$ is cubical, we have $F(u_n'')\supset F(u_n)$. Hence, $z_n\in F(P_1)\setminus P_1$. Without loss of generality we may assume that $z_n\to z\in\cl(F(P_1)\setminus P_1)$. Since
$$
\varrho(z,x')\leq \varrho(z,x)+\varrho(x,x')\leq 2\varepsilon< 1,
$$
we can find $\bar{z}\in\sigma_z\cap\sigma_{x'}\cap\sigma_x$ with $\varrho(x,\bar{z})\leq\varepsilon$. We have $\bar{z}\in\cl(F(P_1)\setminus P_1)\cap P_1$, because $\mathring{\sigma}_z\subset\cl(F(P_1)\setminus P_1)$ and $\sigma_{x'}\subset P_1$. By property (b) of $P$, $\bar{z}\in P_2$. It follows that $x\in\Be(P_2)$, a contradiction.

We shall prove that 
\begin{equation}
\label{eq:inv_part_inclusion_for_vert_enclosure}
\Inv(\Be(N),F^\varepsilon)\subset\Be(\Inv(N,F)).
\end{equation}
Let $x\colon\ZZ\to\Be(N)$ be a solution for $F^\varepsilon$ in $\Be(N)$. For each $x_i\in \Be(N)$, we can choose an $x_i'\in N$ such that $\varrho(x_i,x_i')<\varepsilon$. Since $x_{i+1}\in F^\varepsilon(x_i) = \Be(F(x_i))$, we can take a $z_{i+1}\in F(x_i)$ such that $\varrho(z_{i+1},x_{i+1})<\varepsilon$. We have $\sigma_{z_i}\cap\sigma_{x_i}\cap\sigma_{x'_i}\neq\emptyset$, because $\varrho(z_i,x_i)<\varepsilon$ and $\varrho(x_i,x'_i)<\varepsilon$. Since $F$ has cubical values we get $\sigma_{z_{i+1}}\subset F(x_i)$. For each $i\in\ZZ$ choose a $u_i\in\sigma_{z_i}\cap\sigma_{x_i}\cap\sigma_{x'_i}$. By (\ref{eq:F_eps_prop}) we get $u_{i+1}\in\sigma_{z_{i+1}}\subset F(x_i)\subset F(u_i)$. Since $N$ is cubical and $x'_i\in N$, we get $u_i\in\sigma_{x'_i}\subset N$. Thus, $u_i\in\Inv(N,F)$ and since $\varrho(x_i,u_i)\leq\varrho(x_i,x'_i)\leq\varepsilon$, we get $x_i\in\Be(\Inv(N,F))$. This proves (\ref{eq:inv_part_inclusion_for_vert_enclosure}).

Now, since $\Inv(N,F)$ as an intersection of cubical sets is cubical and $\Inv(N,F)\subset\Int P_1$, we have
$$
\Be(\Inv(N,F))\subset\Int P_1\subset P_1 \subset \Int\Be(P_1).
$$
And, since $\Inv(N,F)\cap P_2 = \emptyset$ and the sets are compact, we have $\Be(\Inv(N,F))\cap\Be(P_2) = \emptyset$. Hence, $\Inv(\Be(N),F^\varepsilon)\subset\Int(\Be(P_1)\setminus \Be(P_2))$, which proves property (c).

In order to prove property (d) it suffices to show that
\begin{equation}
\label{eq:diff_in_neighborhood}
\Be(P_1)\setminus\Be(P_2)\subset N,
\end{equation}
because $N\subset\Int\Be(N)$. Thus, assume that (\ref{eq:diff_in_neighborhood}) is not true and take an $x\in (\Be(P_1)\setminus\Be(P_2))\setminus N$ and choose an $x'\in P_1$ such that $\varrho(x,x')<\varepsilon$. Then $x'\notin P_2$. Let $x''\in\sigma_x\cap\sigma_{x'}$. Since $P$ is cubical, we have $x''\in\sigma_{x'}\subset P_1$. We cannot have $x''\in P_2$, because then $\varrho(x,x'')<\varepsilon$ implies $x\in\Be(P_2)$. Therefore, $x''\in P_1\setminus P_2\subset\Int N$ by property (d) applied to $F$ and $N$. We have $\mathring{\sigma}_x\cap N = \emptyset$, because $x\notin N$ and $N$ is cubical. Thus, $x''\in N\cap \cl\mathring{\sigma}_x\subset\bd N$, a contradiction. 
\qed

For the sake of simplicity in the next theorem for $A\subset X$ we put $A^\varepsilon:=\Be(A)$ and $E(A):=X\setminus\Int A$. 

\begin{thm}
\label{thm:conjugated_index_maps}
Let $F, G\colon X\mto X$ be acyclic upper semicontinuous multivalued maps such that  $F\subset G$. Assume that $N\subset X$ is a cubical isolating neighborhood with respect to $F$, $P$ is a cubical weak index pair in $N$, $N^\varepsilon$ is an isolating neighborhood with respect to $G$, and $P^\varepsilon$ is a weak index pair for $G$ in $N^\varepsilon$. Then the diagram 
$$
\begin{tikzcd}
H^*(P_1,P_2) & \ar{l}[swap]{F^*}H^*(P_1\cup E(N),P_2\cup E(N))\ar{r}{\iota_P^*} &  H^*(P_1,P_2)
\\
 & \ar{u}[swap]{\lambda^*}H^*(P^\varepsilon_1\cup E(N), P^\varepsilon_2\cup E(N)) \ar{d}{\kappa^*} & 
\\
H^*(P^\varepsilon_1,P^\varepsilon_2)\ar{uu}{\alpha^*} & \ar{l}[swap]{G^*} H^*(P_1^\varepsilon\cup E(N^\varepsilon), P_2^\varepsilon\cup E(N^\varepsilon)) \ar{r}{\iota_{P^\varepsilon}^*} & H^*(P^\varepsilon_1,P^\varepsilon_2)\ar{uu}[swap]{\alpha^*}
\end{tikzcd}
$$
commutes and $\alpha^*, \kappa^*, \lambda^*$ are isomorphisms for $0<\varepsilon<\frac{1}{2}$. 
\end{thm}
\proof
Consider the following diagram
$$
\begin{tikzcd}
\arrow[dd,hook,"\alpha"] (P_1,P_2) \arrow[r,"F"]& (P_1\cup E(N),P_2\cup E(N))\arrow[d,hook,"\lambda"] & \arrow[l,hook',swap,"\iota_P"] (P_1,P_2) \arrow[dd,hook,"\alpha"]
\\
 & (P^\varepsilon_1\cup E(N), P^\varepsilon_2\cup E(N)) & 
\\
(P^\varepsilon_1,P^\varepsilon_2) \ar{r}{G} & \arrow[u,hook,swap,"\kappa"] (P_1^\varepsilon\cup E(N^\varepsilon), P_2^\varepsilon\cup E(N^\varepsilon)) & \arrow[l,hook',swap,"\iota_{P^\varepsilon}"](P^\varepsilon_1,P^\varepsilon_2).
\end{tikzcd}
$$
The above diagram commutes up to inclusion, that is $\lambda\circ F\subset \kappa\circ G\circ \alpha$ and $\lambda\circ\iota_P = \kappa\circ\iota_{P^\varepsilon}\circ\alpha$. Inclusions $\iota_P, \iota_{P^\varepsilon}, \kappa$ induce isomorphisms in cohomology by excision. 

Let $\alpha|_{P_i}$ and $\lambda|_{P_i\cup E(N)}$ be restrictions of $\alpha, \lambda$ to appropriate sets, respectively. By Lemma \ref{lem:inclusion_induces_isomorphism}, inclusions $\alpha|_{P_i}\colon P_i \hookrightarrow P_i^\varepsilon$ and $\lambda|_{P_i\cup E(N)}\colon P_i\cup E(N) \hookrightarrow P_i^\varepsilon\cup E(N)$ induce isomorphisms in cohomology for $i=1,2$. Since the following diagram
$$
\begin{tikzcd}
\arrow[d,hook,"\alpha|_{P_2}"] P_2 \arrow[r,hook] & \arrow[d,hook,"\alpha|_{P_1}"] P_1 \arrow[r,hook] & \arrow[d,hook,"\alpha"] (P_1,P_2) 
\\
P_2^\varepsilon \arrow[r,hook] & P_1^\varepsilon  \arrow[r,hook] & (P_1^\varepsilon,P_2^\varepsilon)
\end{tikzcd}
$$
commutes, the diagram
$$
\begin{tikzcd}
... & \ar{l}{} H^q(P_2) & \ar{l}{} H^q(P_1) & \ar{l}{} H^q(P_1,P_2) & \ar{l}{} H^{q-1}(P_2) & \ar{l}{} ... \\
... & \ar{l}{} H^q(P_2^\varepsilon) \ar{u}{(\alpha|_{P_2})^q} & \ar{l}{} H^q(P_1^\varepsilon) \ar{u}{(\alpha|_{P_1})^q} & \ar{l}{} H^q(P_1^\varepsilon,P_2^\varepsilon) \ar{u}{\alpha^q} & \ar{l}{} H^{q-1}(P_2^\varepsilon) \ar{u}{(\alpha|_{P_2})^{q-1}} & \ar{l}{} ... 
\end{tikzcd}
$$
also commutes. By Five Lemma, $\alpha^*$ is an isomorphism. An analogous argument for pairs $(P_1\cup E(N),P_2\cup E(N))$ and $(P^\varepsilon_1\cup E(N), P^\varepsilon_2\cup E(N))$ proves that $\lambda^*$ is an isomorphism too. 
\qed
%-------------------
\begin{thm}
\label{thm:index_maps_f_F_conj}
Let $F\colon X\mto X$ be cubical, upper semicontinuous multivalued map with contractible values. Assume that $N\subset X$ is a cubical isolating neighborhood with respect to $F$, $P$ is a cubical weak index pair in $N$ and $0<\varepsilon<\frac{1}{2}$. Then $a_\varepsilon (F)\neq\emptyset$, and every $\varepsilon$-approximation of $F$ has $N$ as an isolating neighborhood and $R:=\Be(P)\cap N$ as a weak index pair. Moreover, index maps $I_{F_P}$ and $I_{f_{R}}$ are conjugate. 
\end{thm}
\proof
By Propositions \ref{thm:preserving_iso_ne} and \ref{thm:preserving_wip}, $N$ is an isolating neighborhood for $F_\varepsilon$ and $P$ is a weak index pair for $F_\varepsilon$ in $N$. By Lemma \ref{lem:F_sub_epsilon},  $F_\varepsilon$ is upper semicontinuous and has contractible values. Moreover, $F\subset F_\varepsilon$, showing that index maps $I_{F_P}$ and $I_{{F_\varepsilon}_P}$ are conjugate. By Lemma \ref{lem:F_eps_prop} and Theorem \ref{thm:wip_for_vert_enclosure} applied for $F_\varepsilon$ we infer that $\Be(N)$ is an isolating neighborhood for $(F_\varepsilon)^\varepsilon$ and  $\Be(P)$ is a weak index pair for $(F_\varepsilon)^\varepsilon$ in $\Be(N)$. Note that, by Lemma \ref{lem:F_sup_epsilon},  $(F_\varepsilon)^\varepsilon$ is upper semicontinuous and has contractible values. Therefore, Theorem \ref{thm:conjugated_index_maps} applied for maps $F_\varepsilon$ and $(F_\varepsilon)^\varepsilon$, implies that index maps $I_{{F_\varepsilon}_P}$ and $I_{{(F_\varepsilon )^\varepsilon}_{\Be(P)}}$ are conjugate. 

By Proposition \ref{prop:preserving_iso_ne2}, $N$ is an isolating neighborhood for $(F_\varepsilon )^\varepsilon$. Hence, by \cite[Lemma 5.1]{BM2016}, $R$ is a weak index pair for $(F_\varepsilon )^\varepsilon$ in $N$. The diagram
$$
\begin{tikzcd}
\arrow[d,hook,"\alpha"] (R_1,R_2) \arrow[r,"(F_\varepsilon )^\varepsilon"]& (R_1\cup E(N),R_2\cup E(N))\arrow[d,hook,"\id"] & \arrow[l,hook',swap,"\iota_R"] (R_1,R_2) \arrow[d,hook,"\alpha"]\\
(P^\varepsilon_1,P^\varepsilon_2) \ar{r}{(F_\varepsilon )^\varepsilon} & (P_1^\varepsilon\cup E(N), P_2^\varepsilon\cup E(N)) & \arrow[l,hook',swap,"\iota_{P^\varepsilon}"](P^\varepsilon_1,P^\varepsilon_2),
\end{tikzcd}
$$
in which inclusions $\alpha$, $\iota_R$ and $\iota_{P^\varepsilon}$ are excisions, commutes. This, along with the fact that pairs $(R_1,R_2)$ and $(P_1^\varepsilon\cup E(N), P_2^\varepsilon\cup E(N))$ are associate, shows that index maps $I_{{(F_\varepsilon )^\varepsilon}_{\Be(P)}}$ and $I_{{(F_\varepsilon )^\varepsilon}_{R}}$ are conjugate. 

Eventually we infer that  $I_{F_P}$ and $I_{{(F_\varepsilon )^\varepsilon}_{R}}$ are conjugate.

The existence of an $\varepsilon$-approximation $f:X\to X$ of $F$ follows from \cite[Theorem 5.12]{GGK} (cf. e.g. \cite[Theorem 23.8]{G76}). Observe that for an arbitrary $\varepsilon$-approximation $f:X\to X$ of $F$ the inclusion $f\subset(F_\varepsilon)^\varepsilon$ holds. Therefore, index maps $I_{{f_{R}}}$ and $I_{{(F_\varepsilon )^\varepsilon}_{R}}$ are conjugate, and the conclusion follows.
\qed

%----------------------------
%%%%%%%%%%%%%%%%%%%%%%%%%%%%%%%%%%%%%%%%%%%%%%%%%
\section{Index map and its iterates}
\label{sec:indexmap}
Throughout this section we assume that $X$ is a locally compact metrizable space and $f:X\to X$ is a discrete dynamical system.
%-----------------------------------------

For convenience we shall use the notion of associated pairs introduced in \cite{Sr97}. Namely, we say that a pair of paracompact sets $P'=(P' _1, P' _2)$ is {\em associated} to a weak index pair $P$ with respect to $f$, if
\begin{itemize}\label{def:P'_assoc_P}
\item[(a1)] $P\subset P'$;
\item[(a2)] $P_1\setminus P_2= P' _1\setminus P' _2$;
\item[(a3)] $f(P)\subset P'$.
\end{itemize}
Note that if $P'$ is associated to a weak index pair $P$ then the pair of pairs $(P,P')$ is a weak index quadruple in the sense of \cite{Mr06}. Moreover, by (a2) the inclusion $i_{PP'}$ induces an isomorphism in the Alexander-Spanier cohomology, and by (a3), we can consider the restriction $f_{PP'}$ of $f$ to the domain of $P$ as a map of pairs $f_{PP'}:P\to P'$.

Clearly, the pair $T_N(P)$ is associated to $P$. Another pair associated to $P$ is
$$
S_N(P):=(S_{N,1}(P),S_{N,2}(P)):=(P_1\cup (f(P_1)\setminus\Int N),P_2\cup (f(P_1)\setminus\Int N)).
$$
Observe that $S_N(P)$ is the smallest pair associated to $P$, i.e. for any pair $P'$ associated to $P$, we have $S_{N,i}(P)\subset P' _i$. Indeed, for $i=1$ the inclusion follows directly from (a1) and (a3). Note that in order to show the inclusion  $S_{N,2}(P)=P_2\cup (f(P_1)\setminus\Int N)\subset P' _2$ it suffices to verify that $f(P_1)\setminus\Int N\subset P' _2$, as $P_2\subset P' _2$ by (a1). Suppose to the contrary that there exists a $y\in (f(P_1)\setminus\Int N)\setminus P' _2$. Then, by (a3) and (a2), $y\in P' _1\setminus P' _2=P_1\setminus P_2$. However, $P_1\setminus P_2\subset\Int N$; hence $y\in\Int N$, a contradiction.

We have the commutative diagram
$$
\begin{tikzcd}
&(T_{N,1}(P),T_{N,2}(P))&\\[3ex]
(P _1,P _2)\arrow[rd,swap,"f _{PP'}"] \ar{ru}{f _P}\ar{r}{f _{PS(P)}} & (S_{N,1}(P),S_{N,2}(P))\arrow[u, "j_1", swap]\arrow[d, "j_2"] & (P _1,P _2)\arrow[ld,"i_{PP'}"] \ar{l}[swap]{i_{PS(P)}}\ar{lu}[swap]{i_P} \\[3ex]
 & (P'_1,P'_2) &
\end{tikzcd}
$$
in which $i_P$, $i_{PP'}$, $i_{PS(P)}$, $j_1$ and $j_2$ are inclusions. Since any of the pairs in the diagram is associated to $P$, each of the inclusions induces an isomorphism in cohomology. Hence, by the commutativity of the diagram we
obtain $I_{f_P}=f_{PP'}^*\circ (i_{PP'} ^*)^{-1}$. For reference we state this observation as
\begin{prop}\label{prop:I_assoc}
Let $P$ be a weak index pair for $f$ and let $P'$ be a pair associated to $P$. Then
\begin{itemize}
\item[(i)] there is a well defined map of pairs  $f_{PP'}:P\ni x\mapsto f(x)\in P'$,
\item[(ii)] the inclusion $i_{PP'}:P\to P'$ induces an isomorphism in cohomology,
\item[(iii)] $I_{f_P}=f_{PP'}^*\circ (i_{PP'} ^*)^{-1}$.
\end{itemize}
\end{prop}
%---------------------------------------------
\begin{prop}\label{prop:ne_of_N}
Let $M$ be an isolating neighborhood for $f$. For any $n\in\NN$ there exists an open neighborhood $U$ of $\Inv(M,f)$, with $\cl U\subset M$, such that for any $x\in U$ we have
$$
f^k(x)\in\Int M\mbox{ for } k\in I_n.
$$
\end{prop}
\proof
Since $S$ is compact and $f$ is continuous, we can find an open set $U\supset S$, with $\cl U\subset M$, such that $f(U)\cup f^2(U)\cup\dots\cup f^k(U)\subset \Int M$.
\qed\\
%----------------------------------------
The following proposition is straightforward.
\begin{prop}\label{prop:inv_f_subset_inv_fn}
If $N$ is an isolating neighborhood for $f$ then for any $k\in \NN$ we have
\begin{equation}\label{inv_f_subset_inv_fn}
\Inv(N,f)\subset\Inv(N,f^k).
\end{equation}
\end{prop}
Although the converse inclusion is not valid in general, we have the following proposition. 
%----------------------------------------
\begin{prop}\label{prop:iso_ne_f_fn}
Let $S$ be an isolated invariant set with respect to $f$. 
For any $k\in\NN$ there exists an isolating neighborhood $M$ of $S$ such that
\begin{equation}\label{eq:S_inv_f_fn}
\Inv(M,f)=S=\Inv(M,f^k).
\end{equation}
\end{prop}
\proof
Let $\hat{N}$ be an isolating neighborhood of $S$ with respect to $f$. By Proposition \ref{prop:ne_of_N} we can take an open neighborhood $U$ of $S$ such that $\bigcup _{i=1} ^{k}f^k(U)\subset\Int\hat{N}$. Let $M\subset U$ be an isolating neighborhood of $S$. We have $S=\Inv(M,f)\subset \Inv(M,f^k)$. To see the opposite inclusion take an $x\in \Inv(M,f^k)$. Then $f^{ik}(x)\in\Inv(M,f^k)$ for $i\in\ZZ$. But $\Inv(M,f^k)\subset M\subset U$, therefore $f^j(f^{ik}(x))\in\Int \hat{N}$ for $j\in I_k$. Hence, $x\in\Inv(\hat{N},f)=S$.
\qed
%-----------------------------------------------
\begin{prop}\label{prop:wip_f_fn}
Let $S$ be an isolated invariant set for $f$. For any $n\in\NN$ there exist isolating neighborhoods $N\subset M$ of $S$ and weak index pairs $P$ and $Q$, respectively in $N$ and $M$, such that for each $k\in I_n$:
\begin{itemize}
\item[(i)] $P$ is a weak index pair for $S$ and $f^k$,
\item[(ii)] $Q$ is associated to $P$ with respect to $f^k$,
\item[(iii)] $T_N(P)$ is associated to $Q$ with respect to $f$.
\end{itemize}
\end{prop}
\proof Fix an arbitrary $n\in\NN$ and consider an isolating neighborhood $M$ of $S$ satisfying (\ref{eq:S_inv_f_fn}).
Take $U\subset\cl U\subset M$, an open neighborhood of $S$ as in Proposition \ref{prop:ne_of_N}, and a compact set $N\subset U$ with $S\subset\Int N$. Note that such an $N$ is an isolating neighborhood for $f^k$, for each $k\in\{1,2,\dots, n\}$. By \cite[Theorem 4.12]{BM2016} we can find a weak index pair $Q=(Q_1,Q_2)$ for $f$ and $S$ in $M$ such that $Q_1\setminus Q_2\subset\Int N$. Define the pair $P:=(P_1,P_2)$ as the intersection
\begin{equation}\label{eq:def_P}
P:=Q\cap N.
\end{equation}
According to \cite[Lemma 5.1]{BM2016}, $P$ is a weak index pair for $f$ in $N$.
We shall prove that the pairs $P$ and $Q$ satisfy assertions (i), (ii) and (iii).

First we prove that
\begin{equation}\label{eq:fkP_subset_Q}
f^k(P)\subset Q\mbox{ for }k\in I_n.
\end{equation}
We argue by induction with respect to $k$. Since for $i=1,2$ we have $P_i\subset N\subset U$, by Proposition \ref{prop:ne_of_N}, we get $f(P_i)\subset M$. Therefore, $f(P_i)\subset f(Q_i)\cap M\subset Q_i$, as $P_i\subset Q_i$ and $Q_i$ is positively invariant with respect to $f$ and $M$. Next, suppose that for some $k\in I_{n-1}$ we have $f^k(P_i)\subset Q_i$. By Proposition \ref{prop:ne_of_N}, $f^{k+1}(P_i)\subset M$. Consequently, $f^{k+1}(P_i)\subset f(f^k(P_i))\cap M\subset f(Q_i)\cap M\subset Q_i$. This completes the proof of (\ref{eq:fkP_subset_Q}).

We shall prove that $P$ is a weak index pair with respect to each $f^k$, $k\in I_n$. To this end fix an arbitrary $k\in\{2,\dots,n\}$ (recall that for $k=1$ the assertion follows from \cite[Lemma 5.1]{BM2016}). Since $P$ is a weak index pair in $N\subset U$ with respect to $f$, we have $P_1\setminus P_2\subset\Int N$, as well as, $\Inv(N,f^k)=\Inv(N,f)\subset\Int(P_1\setminus P_2)$. This shows that $P$ satisfies properties (c) and (d) of Definition \ref{defn:wip} of a weak index pair for $f^k$. Since property (a) follows easily from (\ref{eq:fkP_subset_Q}), it remains to verify property (b), that is
$
\bd_{f^k}(P_1)\subset P_2.
$
Suppose to the contrary that there exists a $y\in \bd_{f^k}(P_1)\setminus P_2$. Then $y\in P_1\setminus P_2$ and $y\in\cl(f^k(P_1)\setminus P_1)$. Consider a sequence $\{y_n\}\subset f^k(P_1)\setminus P_1$ convergent to $y$. Since $y\in P_1\setminus P_2\subset\Int N$, for sufficiently large $n$ we have $y_n\in\Int N$. Consequently, $y_n\in f^k(P_1)\cap\Int N\subset f^k(P_1)\cap N$, which according to the property (a) of $P$ yields $y_n\in P_1$, a contradiction.

To prove (ii) observe that properties (a1) and (a2) are obvious and (a3) follows from (\ref{eq:fkP_subset_Q}).

We shall show (iii). Since $N\subset M$, by (\ref{eq:def_P}) it follows that $Q\subset T_N(P)$, showing that (a1) is satisfied. Condition (a2) is a direct consequence of (ii) and the fact that $T_N(P)$ is associated to $P$. It remains to verify property (a3). By (\ref{eq:def_P}) and the inclusion $N\subset M$ it follows that $T_M(Q)\subset T_N(P)$. This, along with the obvious inclusion $f(Q)\subset T_M(Q)$, implies $f(Q)\subset T_N(P)$, and completes the proof.
\qed
%---------------------------------------------
\begin{prop}\label{prop:In_f,_I_fn}
Let $N\subset M$ be isolating neighborhoods of $S$. Assume  that $P$ is a weak index pair in $N$ with respect to each $f^k$, $k\in I_p$, and $Q$ is a weak index pair with respect to $f$ in $M$. Moreover, assume  that $Q$ is associated to $P$ with respect to $f^k$, and  $T_N(P)$ is associated to $Q$ with respect to $f$. Then
\begin{equation}\label{eq:In_f,_I_fn}
I_{f^p _P}=I^p_{f _P}.
\end{equation}
\end{prop}
\proof
Fix an arbitrary $k\in I_p$.
Since $Q$ is associated to $P$ with respect to $f^k$, by Proposition \ref{prop:I_assoc} we have
\begin{equation}\label{eq:I_PQ}
I_{f^k _P}=(f^k _{PQ}) ^*\circ (i_{PQ} ^*)^{-1}.
\end{equation}

Note that, for each $k\in I_{p-1}$, we have the commutative diagram
$$
\begin{tikzcd}
(Q_1,Q_2)\ar{dr}{f_{QT_N(P)}}&\ar{l}[swap]{i_{PQ}}(P _1,P _2)\ar{d}{f _P}&
\\[3ex]
(P _1,P _2)\ar{u}{f^k _{PQ}}\ar{r}{f^{k+1} _P} & (T_{N,1}(P),T_{N,2}(P)) & (P _1,P _2)\ar{l}[swap]{i_P}
\end{tikzcd}
$$
in which $i_P$ and $i_{PQ}$ are inclusions. Moreover, the inclusions $i_P$ and $i_{PQ}$ induce isomorphisms in cohomology, as excisions. Therefore, by the commutativity of the diagram and (\ref{eq:I_PQ}) we obtain
$$
\begin{array}{rcl}
I_{f^{k+1} _P}&=&(f^{k+1} _P)^*\circ ({i_P} ^*)^{-1}\\[1ex]
&=&(f^k _{PQ})^*\circ (i_{PQ} ^*)^{-1}\circ (f _{P})^*\circ (i_{P} ^*)^{-1}\\[1ex]
&=&I_{f^{k} _P}\circ I_{f _P}.
\end{array}
$$
Taking into account that the above equality is valid for an arbitrary $k\in I_{p-1}$, the assertion follows by induction.\qed
%--------------------------------------------
\begin{prop}\label{prop:iso_ne_subset}
Assume $N$ is an isolating neighborhood with respect to $f$ and $P$ is a weak index pair for $f$ in $N$. Moreover, assume $N=\bigcup_{i=1}^n N_i$ where $N_i$ are pairwise disjoint compact subsets of $N$. Then, for any $I\subset I_n$, the union $\bigcup _{i\in I} N_i=:N_I$ is an isolating neighborhood for $f$, and $Q:=P\cap N_I$ is a weak index pair for $f$ in $N_I$.
\end{prop}
\proof
Clearly, $N_I$ is compact. Since $\Int N_I=N_I\cap\Int N$ and $N$ is an isolating neighborhood for $f$, we have the inclusions $\Inv(N_I,f)\subset \Inv(N,f)\cap N_I\subset\Int N\cap N_I=\Int N_I$, showing that $N_I$ is an isolating neighborhood for $f$.

We shall verify that $Q$ is a weak index pair in $N_I$. It is obvious that $Q_2\subset Q_1$ are compact subsets of $N_I$. For the proof of condition (a) in Definition \ref{defn:wip} observe that $f(Q_i)\cap N_I\subset f(P_i)\cap N\subset P_i$, hence $f(Q_i)\cap N_I\subset P_i\cap N_I=Q_i$. Moreover, we have the inclusions $\Inv (N_I,f)\subset\Int N_I\cap\Inv(N,f)\subset\Int N_I\cap\Int(P_1\setminus P_2)=\Int(Q_1\setminus Q_2)$, showing that $Q$ satisfies condition (c).
Next, observe that $Q_1\setminus Q_2=(P_1\setminus P_2)\cap N_I\subset\Int N\cap N_I=\Int N_I$, which means that $Q$ satisfies condition (d). We still need to show that $Q$ satisfies property (b). Suppose to the contrary that there exists a $y\in \bd_{f}(Q_1)\setminus Q_2$. Then $y\in Q_1\setminus Q_2$ and $y\in\cl(f(Q_1)\setminus Q_1)$. Thus we can take a sequence $\{y_n\}\subset f(Q_1)\setminus Q_1$ convergent to $y$. By the inclusion $y\in Q_1\setminus Q_2\subset\Int N_I$, it follows that $y_n\in\Int N_I$
for sufficiently large $n$. Consequently, $y_n\in f(Q_1)\cap\Int N_I\subset f(Q_1)\cap N_I$, which according to the positive invariance of $Q_1$ with respect to $f$ and $N_I$ yields $y_n\in Q_1$, a contradiction.
\qed
%----------------------------------------------
\begin{prop}\label{prop:IP_IP1}
Assume that $N$ is an isolating neighborhood  for $f$ and $P$ is a weak index pair in $N$. Moreover, assume $N=N_1\cup N_2$ where $N_1$, $N_2$ are compact disjoint subsets of $N$. Let $P^1:=P\cap N_1$, let $\iota :H^*(P^1)\to H^*(P^1)\times H^*(P^2)$ be the inclusion, and let $\pi:H^*(P^1)\times H^*(P^2)\to H^*(P^1)$ be the projection. Then
\begin{equation}\label{eq:IP_IP1}
I_{f _{P^1}}=\pi\circ I_{f _P}\circ\iota.
\end{equation}
\end{prop}
\proof
By Proposition \ref{prop:iso_ne_subset}, $N_1$ is an isolating neighborhood for $f$, and $P^1$ is a weak index pair in $N^1$. Therefore, we have well defined index maps $I_{f _{P^1}}$ and $I_{f _{P}}$, associated with the weak index pairs $P^1$ and $P$, respectively.

Consider the commutative diagram
$$
\begin{tikzcd}
(P_1,P_2)\ar{r}{f_P}&(T_{N,1}(P),T_{N,2}(P))\ar{d}{k}&(P_1,P_2)\ar{l}[swap]{i_P}
\\[3ex]
(P^1 _1,P^1 _2)\ar{u}{j}\ar{r}{f_{P^1}} & (T_{N_1,1}(P^1),T_{N_1,2}(P^1)) & (P^1 _1,P^1 _2)\ar{l}[swap]{i_{P^1}}\ar{u}[swap]{j}
\end{tikzcd}
$$
in which $i_P$, $i_{P^1}$, $j$ and $k$ are inclusions. Recall that $i_P$ and $i_{P^1}$ induce isomorphisms in cohomology by the strong excision property. By the commutativity of the diagram we obtain
$
j^*\circ f_P ^*\circ(i_P ^*)^{-1}=f_{P^1} ^*\circ(i_{P^1} ^*)^{-1}\circ j^*
$, showing that
\begin{equation}\label{eq:IP_IP1_aux}
j^*\circ I_{f _P}=I_{f _{P^1}}\circ j^*.
\end{equation}
Consider the commutative diagram
$$
\begin{tikzcd}
(P_1,P_2)\ar{rd}{\kappa}&
\\[3ex]
(P^1 _1,P^1 _2)\ar{u}{j}\ar{r}{\lambda}& (P_1,P_2\cup P^2 _1)
\end{tikzcd}
$$
in which $\kappa$ and $\lambda$ are inclusions. Note that $\lambda$ induces an isomorphism in cohomology as an excision. Moreover, by the commutativity of the diagram, $\lambda^*=j^*\circ\kappa ^*$, showing that
\begin{equation}\label{eq:j_star_pi}
j^*\circ(\kappa ^*\circ(\lambda^*)^{-1})=\id _{H^*(P^1 _1, P^1 _2)}.
\end{equation}
Note that $\kappa ^*\circ(\lambda^*)^{-1}:H^*(P^1)\to H^*(P^1)\times H^*(P^2)$ is an inclusion. Thus, $\kappa ^*\circ(\lambda^*)^{-1}=\iota$, and $j^*=\pi$.
Now, (\ref{eq:IP_IP1}) follows from (\ref{eq:IP_IP1_aux}), which completes the proof.
\qed
%%%%%%%%%%%%%%%%%%%%%%%%%%%%%%%%%%%%%%%%%%%%%%%%
\section{Determining orbits via Wa\.zewski property of the Conley index}
Let $X$ be a locally compact metrizable space, and let $f:X\to X$ be a discrete dynamical system. Recall that for $p\geq 2$ we denote by $\ZZ_p:=\{0,1,\dots,p-1\}$ the topological group with the addition modulo $p$ and the  discrete topology. We define the space $\bar{X}:=X\times \ZZ_p$, with the product topology, and dynamical systems $\bar{f},\underline{f}:\bar{X}\to\bar{X}$, by
\begin{equation}\label{eq:def_fbar}
\bar{f}:\bar{X}\ni (x,i)\longmapsto (f(x),i+1)\in\bar{X}
\end{equation}
and
$$
\underline{f}:\bar{X}\ni(x,i)\longmapsto (f(x),i)\in\bar{X},
$$
respectively. Consider the homeomorphism
$$
l:\bar{X}\ni(x,i)\longmapsto (x,i+1)\in\bar{X}
$$
and observe that we have
\begin{equation}\label{eq:bar_f_und_f}
\bar{f}=\underline{f}\circ l=l\circ\underline{f}.
\end{equation}
Given $A\subset X$, by $\bar{A}$ we shall denote the set $A\times \ZZ_p$.
\begin{prop}\label{prop:barN_iso_ne}
If $N$ is an isolating neighborhood for $f$ then $\bar{N}$ is an isolating neighborhood for both $\bar{f}$ and \underline{f}.
Moreover, if $P$ is a weak index pair for $f$ in $N$ then $\bar{P}$ is a weak index pair in $\bar{N}$ for both $\bar{f}$ and \underline{f}.
\end{prop}
\proof
Consider the dynamical system $\bar{f}$. Clearly $\bar{N}$ is compact. We shall verify that $\Inv(\bar{N},\bar{f})\subset\Int\bar{N}$. To this end consider $\bar{x}=(x,i)\in \Inv(\bar{N},\bar{f})$. Let $\bar{\sigma}:\ZZ\to\bar{X}$ be a solution for $\bar{f}$ passing through $\bar{x}$, which is contained in $\bar{N}$, that is $\bar{\sigma}(0)=\bar{x}$, $\bar{\sigma}(\ZZ)\subset\bar{N}$ and $\bar{\sigma}(k+1)=\bar{f}(\bar{\sigma}(k))$ for $k\in\ZZ$. Define $\sigma:\ZZ\to X$ by $\sigma(k):=p(\bar{\sigma}(k))$ for $k\in\ZZ$, where $p:\bar{X}\ni(x,i)\mapsto x\in X$ denotes the projection. One can easily see that $\sigma$ is a solution for $f$ through $x$ in $N$. Therefore, $x\in\Int N$, as $N$ is an isolating neighborhood for $f$. This shows that $\bar{x}=(x,i)\in\Int N\times \ZZ_p=\Int\bar{N}$, and completes the proof.

The verification that $\bar{P}$ is a weak index pair for $\bar{f}$ and $\bar{N}$ is straightforward.

The proof for $\underline{f}$ is similar.
\qed\\

%----------------------------------------------
For $i\in\ZZ_p$ define the map
\begin{equation}\label{eq:def_mu}
\mu_i:X\ni x\mapsto (x,i)\in X\times\{i\}.
\end{equation}
The following proposition is straightforward.
\begin{prop}\label{lem:f_underline_i}
Assume that $N$ is an isolating neighborhood for $f$, and $P$ is a weak index pair in $N$. For any $i\in\ZZ_p$ the set $N\times\{i\}$ is an isolating neighborhood for $\underline{f}$, and $P\times\{i\}$ is a weak index pair in $N\times\{i\}$. Moreover,
$$
I_{f_P}\circ\mu_i ^*=\mu_i ^* \circ I_{\underline{f}_{P\times\{i\}}}.
$$
\end{prop}
%---------------------------------------------
\begin{prop}\label{prop:I_bar}
Assume that $N$ is an isolating neighborhood for $f$, and $P$ is a weak index pair in $N$. We have
$$
(\times _{i=0} ^{p-1}\mu _i ^*)\circ I_{\bar{f}_{\bar{P}}}=(\times _{i=0} ^{p-1} I_{f_P})\circ (\times _{i=0} ^{p-1}\mu _{i+1} ^*),
$$
where $\times _{i=0} ^{p-1}\mu _i ^*:H^*(\bar{X})\to \bigtimes _{i=0} ^{p-1}H^*(X)$.
\end{prop}
\proof
By Proposition \ref{prop:barN_iso_ne} the pair $\bar{P}$ is a weak index pair with respect to $\underline{f}$. Therefore, the restriction $\underline{f}_{\bar P}$ of $\underline{f}$ to the domain $\bar{P}$ is a map of pairs
$$
\underline{f}_{\bar{P}}:\bar{P}\to T_{\bar{N}}(\bar{P}).
$$
We claim that
\begin{equation}\label{eq:p63e}
I_{\bar{f}_{\bar{P}}}=I_{\underline{f}_{\bar{P}}}\circ l^*.
\end{equation}
Indeed, note that $l\circ i_{\bar{P}}=i_{\bar{P}}\circ l$. Hence, $l^*\circ(i_{\bar{P}}^*)^{-1}=(i_{\bar{P}}^*)^{-1}\circ l^*$ and, by the second equality in (\ref{eq:bar_f_und_f}), we get
$$
\begin{array}{rcl}
I_{\bar{f}_{\bar{P}}}&=&\bar{f}_{\bar{P}}^*\circ (i_{\bar{P}}^*)^{-1}\\[1ex]
&=&(l\circ\underline{f}_{\bar{P}})^*\circ (i_{\bar{P}}^*)^{-1}\\[1ex]
&=&\underline{f}_{\bar{P}}^*\circ l^*\circ (i_{\bar{P}}^*)^{-1}\\[1ex] 
&=&\underline{f}_{\bar{P}}^*\circ (i_{\bar{P}}^*)^{-1}\circ l^*\\[1ex]
&=&I_{\underline{f}_{\bar{P}}}\circ l^*.
\end{array}
$$
For $i\in\ZZ_p$ denote by $l_i$ the restriction of $l$ to the domain $X\times\{i\}$, and observe that $l_i\circ\mu_i=\mu_{i+1}$. Hence, $\mu_i ^*\circ l_i ^*=\mu_{i+1} ^*$ and we have
\begin{equation}\label{eq:mui_l_mui+1}
(\times _{i=0} ^{p-1}\mu _i ^*)\circ l^*=(\times _{i=0} ^{p-1}\mu _i ^*)\circ(\times _{i=0} ^{p-1}l_i ^*)=\times _{i=0} ^{p-1}(\mu _i ^*\circ l_i ^*)=\times _{i=0} ^{p-1}\mu _{i+1} ^*.
\end{equation}
Therefore, according to (\ref{eq:p63e}), in order to complete the proof it suffices to verify that
$$
(\times _{i=0} ^{p-1}\mu _i ^*)\circ I_{\underline{f}_{\bar{P}}}=(\times _{i=0} ^{p-1}I_{f_P})\circ (\times _{i=0} ^{p-1}\mu _i ^*).
$$
Since $\bar{P}$ is a union of pairwise disjoint sets $P\times\{i\}$, we have the product decomposition of $H^*(\bar{P})=\bigtimes _{i=0} ^{p-1} H^*(P\times\{i\})$. Similarly,
$H^*(T_{\bar{N}}(\bar{P}))=\bigtimes _{i=0} ^{p-1}H^*(T_N(P)\times \{i\})$, as the sets
$T_N(P)\times \{i\}$ are pairwise disjoint.
According to the definition of $\underline{f}$ and Proposition \ref{lem:f_underline_i}, we can consider the restriction $\underline{f}_{P\times \{i\}}$ of $\underline{f}$ to the domain $P\times \{i\}$ as a map of pairs
$$
\underline{f}_{P\times \{i\}}:P\times \{i\}\to T_N(P)\times \{i\}.
$$
Thus, we have
$$
\underline{f}_{\bar{P}} ^*=\times _{i=0} ^{p-1}\underline{f}_{P\times\{i\}} ^*.
$$
Similarly,
$$
i_{\bar{P}} ^*= \times _{i=0} ^{p-1}i_{P\times\{i\}} ^*.
$$
Consequently,
$$
\begin{array}{rcl}
I_{\underline{f}_{\bar{P}}}&=&\underline{f}_{\bar{P}} ^*\circ\left(i_{\bar{P}}^*\right)^{-1}\\[2ex]
&=&\left(\times _{i=0} ^{p-1}\underline{f}_{P\times\{i\}} ^*\right)\circ\left(\times _{i=0} ^{p-1}\left(i_{P\times\{i\}} ^*\right)^{-1}\right)\\[2ex]
&=&\times _{i=0} ^{p-1}\left(\underline{f}_{P\times\{i\}} ^*\circ\left(i_{P\times\{i\}}^*\right)^{-1}\right)\\[2ex]
&=&\times _{i=0} ^{p-1}I_{\underline{f}_{P\times\{i\}}}.
\end{array}
$$
Now, by Proposition \ref{lem:f_underline_i}, we obtain
$$
\begin{array}{rcl}
\left(\times _{i=0} ^{p-1}\mu _i ^*\right)\circ I_{\underline{f}_{\bar{P}}}&=&\left(\times _{i=0} ^{p-1}\mu _i ^*\right)\circ \left(\times _{i=0} ^{p-1}I_{\underline{f}_{P\times\{i\}}}\right)\\[2ex]
&=&\times _{i=0} ^{p-1}\left(\mu _i ^*\circ I_{\underline{f}_{P\times\{i\}}}\right)\\[2ex]
&=&\times _{i=0} ^{p-1}\left(I_{{f}_{P}}\circ \mu _i ^*\right)\\[2ex]
&=&\left(\times _{i=0} ^{p-1}I_{f_P}\right)\circ \left(\times _{i=0} ^{p-1}\mu _i ^*\right),
\end{array}
$$
which completes the proof.\qed\\

From now on we assume that $N=\bigcup_{i=1}^n N_i$, where $N_i$ are pairwise disjoint compact subsets of $N$, is an isolating neighborhood with respect to $f$, and $P$ is a weak index pair for $f$ in $N$. Denote $P^i:=P\cap N_i$. 
Let $p\in\NN$ and let $\sigma:=(\sigma_0,\dots,\sigma_{p-1})\in I_n ^{\ZZ_p}$. Consider endomorphism $I_\sigma:\bigtimes _{i=0} ^{p-1}H^*(P^{\sigma_i})\to \bigtimes _{i=0} ^{p-1}H^*(P^{\sigma_i})$ given by
\begin{equation}\label{eq:Isigma}
I_{\sigma}:=\times _{i=0} ^{p-1}\left(\pi _{\sigma _i}\circ I_{f_P}\circ\iota_{\sigma_{i+1}}\right),
\end{equation}
where $\pi _i: H^*(P)\to H^*(P^{i})$ are projections, and $\iota _i:H^*(P^{i})\to H^*(P)$ are inclusions. 

Consider the dynamical system $\bar{f}$ on $\bar{X}$ given by (\ref{eq:def_fbar}). For $\sigma\in I_n ^{\ZZ_p}$ set 
$$
N_{\sigma}:=\bigcup _{i=0} ^{p-1}(N_{\sigma_i}\times\{i\})
$$ 
and let
\begin{equation}\label{eq:def_S_sigma}
S_\sigma:=\Inv(N_\sigma, \bar{f}).
\end{equation}
%------------------------------------------
\begin{prop}\label{prop:Isigma_IR_conjugate}
The set $S_\sigma$ is an isolated invariant set for $\bar{f}$, $N_\sigma$ is its isolating neighborhood, and there exists a weak index pair $R$ for $\bar{f}$ and $S_\sigma$ such that 
\begin{equation}\label{eq:prop_Is_Ir}
I_\sigma\circ\left(\times _{i=0} ^{p-1}\mu _{i+1} ^*\right)=\left(\times _{i=0} ^{p-1}\mu _{i} ^*\right)\circ I_{\bar{f}_R}.
\end{equation}
Moreover, $I^p _{\bar{f}_R}$ and $I^p _\sigma$ are conjugate.
\end{prop}
\proof
First note that, by Proposition \ref{prop:barN_iso_ne}, $\bar{N}$ is an isolating neighborhood for $\bar{f}$, and $\bar{P}$ is a weak index pair in $\bar{N}$. Clearly, $N_\sigma$ is a compact subset of $\bar{N}$; hence, according to Proposition \ref{prop:iso_ne_subset}, $N_\sigma$ is an isolating neighborhood for $\bar{f}$, and $R:=\bar{P}\cap N_\sigma$ is a weak index pair in $N_\sigma$. Therefore, we have a well defined index map $I_{\bar{f}_R}$ for $\bar{f}$, associated with the weak index pair $R$.

We shall prove that $I_{\bar{f}_R}$ and $I_\sigma$ satisfy (\ref{eq:prop_Is_Ir}). To this end consider projections $\bar{\pi}_{k,i}:H^*(P\times\{i\})\to H^*(P^k\times\{i\})$ and the inclusions $\bar{\iota}_{k,i}:H^*(P^k\times\{i\})\to H^*(P\times\{i\})$, for $k\in I_n$ and $i\in\ZZ_p$. One can observe that, for any $i\in\ZZ_p$, we have
$$
\mu _i ^*\circ \bar{\pi}_{k,i}=\pi_{k}\circ \mu _i ^*
$$
and 
$$
\mu _i ^*\circ \bar{\iota}_{k,i}=\iota_{k}\circ \mu _i ^*.
$$
Using the above identities and Proposition \ref{prop:I_bar} we obtain
$$
\begin{array}{rcl}
I_\sigma\circ\left(\times _{i=0} ^{p-1}\mu _{i+1} ^*\right)&=&\left(\times _{i=0} ^{p-1}\pi_{\sigma_i}\right)\circ \left(\times _{i=0} ^{p-1}I_{f_P}\right)\circ\left(\times _{i=0} ^{p-1}\iota_{\sigma_{i+1}}\right)\circ\left(\times _{i=0} ^{p-1}\mu _{i+1} ^*\right)\\[1ex]
&=&\left(\times _{i=0} ^{p-1}\pi_{\sigma_i}\right)\circ \left(\times _{i=0} ^{p-1}I_{f_P}\right)\circ\left(\times _{i=0} ^{p-1}\left(\iota_{\sigma_{i+1}}\circ\mu _{i+1} ^*\right)\right)\\[1ex]
&=&\left(\times _{i=0} ^{p-1}\pi_{\sigma_i}\right)\circ \left(\times _{i=0} ^{p-1}I_{f_P}\right)\circ\left(\times _{i=0} ^{p-1}\left(\mu _{i+1} ^*\circ\bar{\iota}_{\sigma_{i+1},i+1}\right)\right)\\[1ex]
&=&\left(\times _{i=0} ^{p-1}\pi_{\sigma_i}\right)\circ \left(\times _{i=0} ^{p-1}I_{f_P}\right)\circ\left(\times _{i=0} ^{p-1}\mu _{i+1} ^*\right)\circ\left(\times _{i=0} ^{p-1}\bar{\iota}_{\sigma_{i+1},i+1}\right)\\[1ex]
&=&\left(\times _{i=0} ^{p-1}\pi_{\sigma_i}\right)\circ\left(\times _{i=0} ^{p-1}\mu _{i} ^*\right)\circ I_{\bar{f}_{\bar{P}}}\circ\left(\times _{i=0} ^{p-1}\bar{\iota}_{\sigma_{i+1},i+1}\right)\\[1ex]
&=&\left(\times _{i=0} ^{p-1}\left(\pi_{\sigma_i}\circ\mu _{i} ^*\right)\right)\circ I_{\bar{f}_{\bar{P}}}\circ\left(\times _{i=0} ^{p-1}\bar{\iota}_{\sigma_{i+1},i+1}\right)\\[1ex]
&=&\left(\times _{i=0} ^{p-1}\left(\mu _{i} ^*\circ\bar{\pi}_{\sigma_i,i}\right)\right)\circ I_{\bar{f}_{\bar{P}}}\circ\left(\times _{i=0} ^{p-1}\bar{\iota}_{\sigma_{i+1},i+1}\right)\\[1ex]
&=&\left(\times _{i=0} ^{p-1}\mu _{i} ^*\right)\circ\left(\times _{i=0} ^{p-1}\bar{\pi}_{\sigma_i,i}\right)\circ I_{\bar{f}_{\bar{P}}}\circ\left(\times _{i=0} ^{p-1}\bar{\iota}_{\sigma_{i+1},i+1}\right).
\end{array}
$$
Note that $\times _{i=0} ^{p-1}\bar{\pi}_{\sigma_i,i}$ is the projection of $H^*(\bar{P})$ onto $H^*(R)$, and $\times _{i=0} ^{p-1}\bar{\iota}_{\sigma_{i+1},i+1}$ is the inclusion of $H^*(R)$ into $H^*(\bar{P})$. Hence, applying 
Proposition \ref{prop:IP_IP1} we get (\ref{eq:prop_Is_Ir}).

We shall prove that
\begin{equation}\label{eq:conj_Is}
I^p_\sigma\circ\left(\times _{i=0} ^{p-1}\mu _{i+1} ^*\right)=\left(\times _{i=0} ^{p-1}\mu _{i+1} ^*\right)\circ I^p _{\bar{f}_R}.
\end{equation}
Note that $\bar{f}_R\circ l=l\circ\bar{f}_R$. Hence,  $(l^*)^{-1}\circ \bar{f}^* _R=\bar{f}^* _R\circ (l^*)^{-1}$. Similarly, $(l^*)^{-1}\circ (i^* _R)^{-1}=(i^* _R)^{-1}\circ (l^*)^{-1}$, as $i_R\circ l=l\circ i_R$. We have
$$
\begin{array}{rcl}
(l^*)^{-1}\circ I_{\bar{f}_R}&=&(l^*)^{-1}\circ\bar{f}^* _R\circ (i^* _R)^{-1}\\[1ex]
&=&\bar{f}^* _R\circ (l^*)^{-1}\circ (i^* _R)^{-1}\\[1ex]
&=&\bar{f}^* _R\circ (i^* _R)^{-1}\circ (l^*)^{-1}\\[1ex]
&=&I_{\bar{f}_R}\circ(l^*)^{-1}.
\end{array}
$$
Therefore, using (\ref{eq:prop_Is_Ir}) and (\ref{eq:mui_l_mui+1}), we obtain
$$
\begin{array}{rcl}
I^p _\sigma\circ\left(\times _{i=0} ^{p-1}\mu _{i+1} ^*\right)&=&I^{p-1} _\sigma\circ\left(\times _{i=0} ^{p-1}\mu _{i} ^*\right)\circ I_{\bar{f}_R}\\[1ex]
&=&I^{p-1} _\sigma\circ\left(\times _{i=0} ^{p-1}\mu _{i+1} ^*\right)\circ (l^*)^{-1}\circ I_{\bar{f}_R}\\[1ex]
&=&I^{p-1} _\sigma\circ\left(\times _{i=0} ^{p-1}\mu _{i+1} ^*\right)\circ I_{\bar{f}_R}\circ (l^*)^{-1}.
\end{array}
$$
Now, by the reverse induction with respect to $p$, and the fact that $\left((l^*)^{-1}\right)^p$ is the identity, we get (\ref{eq:conj_Is}).
This shows that $I^p _{\bar{f}_R}$ and $I^p _\sigma$ are conjugate, and completes the proof.
\qed

%------------------------------------------
We are ready to present main theorems of this section. They show that from the index map for $f$, itself, we can extract an information which is sufficient to justify the existence of an orbit of $f$, passing through the components of $N$ in a given order.
\begin{thm}\label{thm:Wazewski}
Assume that $N=\bigcup_{i=1}^n N_i$, where $N_i$ are pairwise disjoint compact subsets of $N$, is an isolating neighborhood with respect to $f$, and $P$ is a weak index pair for $f$ in $N$. Let $p\in\NN$ and let $\sigma:=(\sigma_0,\dots,\sigma_{p-1})\in I_n ^{\ZZ_p}$. If the endomorphism $I_\sigma$ given by (\ref{eq:Isigma}) is not nilpotent then there exists  
a trajectory $\tau:\ZZ\to \Inv(\bigcup _{i=0} ^{p-1}N_{\sigma_i},f)$ for $f$, such that $\tau(i+kp)\in N_{\sigma _i}$, for $i\in I_p$, $k\in\ZZ$.
\end{thm}
\proof
By Proposition \ref{prop:Isigma_IR_conjugate}, $S_\sigma =\Inv(N_\sigma, \bar{f})$ is an isolated invariant set for $\bar{f}$. Thus, we have a well-defined Conley index $C(S_\sigma,\bar{f})$ for $S_\sigma$ and $\bar{f}$. Note that, by Proposition \ref{prop:Isigma_IR_conjugate}, 
there exists a weak index pair $R$ in $\bar{X}$ for $\bar{f}$ and $S_\sigma$, such that $I_{\bar{f}_R}$ and $I_\sigma$ satisfy (\ref{eq:prop_Is_Ir}). Since $I_\sigma$ is not nilpotent, then so is $I_{\bar{f}_R}$. Consequently, $C(S_\sigma,\bar{f})\neq 0$. By the Wa\.zewski property of the Conley index (cf. \cite[Proposition 2.10]{M90}), it follows that $S_\sigma\neq\emptyset$. According to definition (\ref{eq:def_fbar}) of $\bar{f}$, there exists an $x\in N_{\sigma_0}$ such that $(x,0)\in S_\sigma$. Let $\eta:\ZZ\to S_\sigma$ be a trajectory for $\bar{f}$ in $S_\sigma$ through $(x,0)$. One easily verifies that then $\tau:=p\circ\eta$, where $p:\bar{X}\ni (x,i)\mapsto x\in X$ is the projection, is a trajectory for $f$ satisfying the assertion. 
\qed

%%%%%%%%%%%%%%%%%%%%%%%%%%%%%%%%%%%%%%%%%%%%%%%%
For a given $i\in I_n$ define endomorphism $g_i:H^*(P)\to H^*(P)$, by
\begin{equation}\label{eq:def_g}
g_i:=I_{f_P}\circ \iota_i\circ\pi_i.
\end{equation}
We are going to prove the following theorem which may be viewed as a counterpart of Theorem \ref{thm:Wazewski} expressed in terms of compositions of endomorphisms $g_i$.
%--------
\begin{thm}\label{thm:Wazewski_g}
Assume that $N=\bigcup_{i=1}^n N_i$, where $N_i$ are pairwise disjoint compact subsets of $N$, is an isolating neighborhood with respect to $f$, and $P$ is a weak index pair for $f$ in $N$. Let $p\in\NN$, let $\sigma:=(\sigma_0,\dots,\sigma_{p-1})\in I_n ^{\ZZ_p}$ and let endomorphisms $g_i$ be given by (\ref{eq:def_g}). If the composition $g_{\sigma_0}\circ\dots\circ g_{\sigma _{p-1}}$ is not nilpotent then there exists 
a trajectory $\tau:\ZZ\to \Inv(\bigcup _{i=0} ^{p-1}N_{\sigma_i},f)$ for $f$, such that $\tau(i+kp)\in N_{\sigma _i}$, for $i\in I_p$, $k\in\ZZ$.
\end{thm}
For its proof we need an auxiliary lemma. Consider the projections 
$$
r_i:\bigtimes _{i=0} ^{p-1}H^*(P^{\sigma_i})\to H^*(P^{i})
$$ 
and the inclusions
$$
m_i:H^*(P^{i})\to \bigtimes _{i=0} ^{p-1}H^*(P^{\sigma_i}).
$$
Let $h_i:\bigtimes _{i=0} ^{p-1}H^*(P^{\sigma_i})\to \bigtimes _{i=0} ^{p-1}H^*(P^{\sigma_i})$ be given by
\begin{equation}\label{eq:def_h}
h_i:=I_\sigma\circ m_i\circ r_i.
\end{equation}
Let $\perm (\ZZ_p)$ and $\cycle (\ZZ_p)\subset \perm (\ZZ_p)$ stand for the sets of all permutations and all cyclic translations of $\ZZ_p$, respectively.
%-----------------------------------------
\begin{lem}\label{lem:Ip_h_g}
Assume $I_\sigma$, $g_i$ and $h_i$ are given by (\ref{eq:Isigma}), (\ref{eq:def_g}) and (\ref{eq:def_h}), respectively. Then:
\begin{itemize}
\item[(i)] $I_\sigma ^p=\Sigma _{s\in \cycle (\ZZ_p)}(h_{\sigma_{s(0)}}\circ\dots\circ h_{\sigma_{s(p-1)}})$,
\item[(ii)] $
h_{\sigma_0}\circ\dots\circ h_{\sigma_{p-1}}=m_{\sigma_{p-1}}\circ\pi_{\sigma_{p-1}}\circ g_{\sigma_{0}}\circ\dots\circ g_{\sigma_{p-1}}\circ\iota_{\sigma_{p-1}}\circ r_{\sigma_{p-1}}
$.
\end{itemize}
\end{lem}
\proof
One can observe that 
\begin{equation}\label{eq:Isigma_sum_hi}
I_\sigma=\Sigma _{i=0} ^{p-1}h_{\sigma_i}.
\end{equation}
Since $h_{\sigma_j}\circ h_{\sigma_i}=0$ whenever $i-j\neq 1$, $i,j\in\ZZ_p$, by (\ref{eq:Isigma_sum_hi}) we have 
$$
\begin{array}{rcl}
I_\sigma ^p&=&\Sigma _{s\in \perm (\ZZ_p)}\left(h_{\sigma_{s(0)}}\circ\dots\circ h_{\sigma_{s(p-1)}}\right)\\[1ex]
&=&\Sigma _{s\in \cycle (\ZZ_p)}\left(h_{\sigma_{s(0)}}\circ\dots\circ h_{\sigma_{s(p-1)}}\right),
\end{array}
$$
which completes the proof of (i).

For the proof of (ii) first observe that, according to the definitions (\ref{eq:Isigma}) and (\ref{eq:def_g}) of $I_\sigma$  and $g_{i}$, respectively, we have the following representation of endomorphisms $h_i$ given by (\ref{eq:def_h})
\begin{equation}\label{eq:expr_h_by_g}
h_{\sigma_{i+1}}=m_{\sigma_i}\circ\pi_{\sigma_i}\circ g_{\sigma_{i+1}}\circ \iota_{\sigma_{i+1}}\circ r_{\sigma_{i+1}}.
\end{equation}
It is straightforward to see that, for each $i\in I_n$, we have
\begin{equation}\label{eq:g_comp=g}
g_{\sigma_i}\circ\iota_{\sigma_i}\circ r_{\sigma_i}\circ m_{\sigma_i}\circ\pi_{\sigma_i}=g_{\sigma_i}.
\end{equation}
Therefore, using (\ref{eq:expr_h_by_g}), we obtain 
$$
h_{\sigma_0}\circ\dots\circ h_{\sigma_{p-1}}=m_{\sigma_{p-1}}\circ\pi_{\sigma_{p-1}}\circ g_{\sigma_{0}}\circ\dots\circ g_{\sigma_{p-1}}\circ\iota_{\sigma_{p-1}}\circ r_{\sigma_{p-1}}.
$$
This completes the proof.
\qed\\

{\bf Proof of Theorem \ref{thm:Wazewski_g}.} According to Theorem \ref{thm:Wazewski} it suffices to show that $I_\sigma ^p$ is not nilpotent. For contradiction suppose that $I_\sigma ^p$ is nilpotent and consider $k\in\NN$ such that $I_\sigma ^{kp}=0$. Note that, by Lemma \ref{lem:Ip_h_g} and the fact that $h_{\sigma_j}\circ h_{\sigma_i}=0$ for $i-j\neq 1$, $i,j\in\ZZ_p$, it follows that
$$
\begin{array}{rcl}
I_\sigma ^{kp}&=&\left(\Sigma _{s\in \cycle (\ZZ_p)}\left(h_{\sigma_{s(0)}}\circ\dots\circ h_{\sigma_{s(p-1)}}\right)\right)^k\\[1ex]
&=&\Sigma _{s\in \cycle (\ZZ_p)}\left(h_{\sigma_{s(0)}}\circ\dots\circ h_{\sigma_{s(p-1)}}\right)^k.
\end{array}
$$
Hence, according to definition (\ref{eq:def_h}) of $h_i$, for each $s\in \cycle (\ZZ_p)$ we have  $(h_{\sigma_{s(0)}}\circ\dots\circ h_{\sigma_{s(p-1)}})^k=0$. In particular, $(h_{\sigma_{0}}\circ\dots\circ h_{\sigma_{p-1}})^k=0$. Consequently, by Lemma \ref{lem:Ip_h_g}(ii) and (\ref{eq:g_comp=g}), we obtain
$$
m_{\sigma_{p-1}}\circ\pi_{\sigma_{p-1}}\circ (g_{\sigma_{0}}\circ\dots\circ g_{\sigma_{p-1}})^k\circ\iota_{\sigma_{p-1}}\circ r_{\sigma_{p-1}}=0,
$$
which implies $(g_{\sigma_{0}}\circ\dots\circ g_{\sigma_{p-1}})^k=0$, a contradiction. This completes the proof.\qed. 

%%%%%%%%%%%%%%%%%%%%%%%%%%%%%%%%%%%%%%%%%%%%%%%%
\section{Determining periodic orbits via Lefschetz-type fixed point theorem}
\label{sec:lefschetz}
We will continue to deal with determining orbits passing through the disjoint components of an isolating neighborhood in a prescribed fashion. Now we focus our attention on periodic orbits.
  
Throughout this section we use the notation introduced in the preceding section. 

Let $\varphi=\{\varphi _i\}$ be an endomorphism of degree zero of a graded vector space $V=\{V_i\}$. Recall that $\varphi$ is called a {\em Leray endomorphism} provided the quotient space $V':=V/N(\varphi)$, where $N(\varphi):=\bigcup\{\varphi^{-n}(0)\;|\;n=1,2,\dots\}$, is of a finite type. For such a $\varphi$ we define its trace as a trace of an induced endomorphism $\varphi':V'\to V'$, i.e. $\tr(\varphi):=\tr(\varphi')$, and the {\em (generalized) Lefschetz number}, by
$$
\Lambda(\varphi):=\sum _{i=0} ^\infty (-1)^i\tr(\varphi_i).
$$ 
It is worth to mention the case of endomorphisms $\varphi,\psi$ of graded vector spaces $V$ and $W$, respectively, such that $\varphi=hg$ and $\psi=gh$ for some morphisms $g:V\to W$ and $h:W\to V$. If one of such endomorphisms is a Leray endomorphism then so is the other, and $\Lambda(\varphi^k)=\Lambda(\psi^k)$ for all $k\in\NN$ (cf. \cite{Gr72}, \cite[Proposition 2]{Mr90}). It applies, in particular, if $\varphi$ and $\psi$ are {\em conjugate}, that is, there exists an isomorphism $g:V\to W$ such that $g\varphi=\psi g$.

%-----------------------------------------
The following Proposition shows that the Lefschetz number of an index map is independent of the choice of a weak index pair.
\begin{prop}\label{prop:Lef_independent_of_P}
Let $S$ be an isolated invariant set for $f$ and let $P$ and $R$ be arbitrary weak index pairs for $f$ and $S$. Then, for every $k\in\NN$, if $\Lambda(I^k _{f _P})$ is well defined, then so is $\Lambda(I^k _{f_R})$ and we have
\begin{equation}\label{eq:Lef_independent_of_P}
\Lambda(I^k _{f _P})=\Lambda(I^k _{f_R}).
\end{equation} 
\end{prop}
\proof
By \cite[Theorem 6.4]{BM2016} and its proof it follows that there exists a sequence $I_{f_P}=I_1,I_2,\dots,I_k=I_{f_R}$ of endomorphisms, with the property that each two consecutive endomorphisms, $I_i$ and $I_{i+1}$, are linked in the sense of \cite[Proposition 2]{Mr90}. Hence, the assertion follows.
\qed
%------------------------------------------
\begin{prop}\label{thm:periodic_trajectory1}
For any weak index pair $Q$ for $\bar{f}^p$ and $S_\sigma$ given by (\ref{eq:def_S_sigma}), if $\Lambda(I_{\bar{f}^p _{Q}})$ is well defined then so is $\Lambda(I_\sigma ^p)$ and we have
$$
\Lambda(I_{\bar{f}^p _{Q}})=\Lambda(I_\sigma ^p).
$$
\end{prop}
\proof
By Proposition \ref{prop:iso_ne_f_fn}, $S_\sigma$ is an isolated invariant set with respect to both $\bar{f}$ and $\bar{f}^p$. Moreover, 
according to Proposition \ref{prop:wip_f_fn}, we can take a pair  $P'$, which is a weak index pair for each $\bar{f}^k$, $k\in I_p$ and $S_\sigma$, and satisfies all the assumptions of Proposition \ref{prop:In_f,_I_fn}.
Then Proposition \ref{prop:In_f,_I_fn} implies that
\begin{equation}\label{eq:Lef_I_fP'p=Lef_Ip_fP'}
\Lambda(I_{\bar{f}^p _{P'}})=\Lambda(I^p_{\bar{f} _{P'}}).
\end{equation}
Since $Q$ is a weak index pair for $\bar{f}^p$ and $S_\sigma$, and so is $P'$, by Proposition \ref{prop:Lef_independent_of_P}, we get  
\begin{equation}\label{eq:Lef_IP'_fp=Lef_IQ_fp}
\Lambda(I_{\bar{f}^p _{Q}})=\Lambda(I_{\bar{f}^p _{P'}}).
\end{equation}
According to Proposition \ref{prop:Isigma_IR_conjugate} we can take a weak index pair $R$ for $\bar{f}$ and $S_\sigma$, such that $I^p _\sigma$ and $I^p _{\bar{f}_R}$ are conjugate; hence,
\begin{equation}\label{eq:IR^p_Isigma^p}
\Lambda(I_{\bar{f} _{R}} ^p)=\Lambda(I_\sigma ^p).
\end{equation}
Note that both $P'$ and $R$ are weak index pairs for $\bar{f}$ and $S_\sigma$. Therefore, applying  Proposition \ref{prop:Lef_independent_of_P} once again, we have
\begin{equation}\label{eq:Lef_IpR=Lef_IpP'}
\Lambda(I_{\bar{f} _{P'}} ^p)=\Lambda(I_{\bar{f} _{R}} ^p).
\end{equation}
Now, the assertion follows from (\ref{eq:Lef_IP'_fp=Lef_IQ_fp}), (\ref{eq:Lef_I_fP'p=Lef_Ip_fP'}),  (\ref{eq:Lef_IpR=Lef_IpP'}), and
(\ref{eq:IR^p_Isigma^p}).   
\qed\\

%------------------------------------------
Note that $\bar{f}^p$ maps $X\times\{i\}\subset \bar{X}$ into itself, for any $i\in I_p$. Therefore, the following proposition is straightforward.
\begin{prop}\label{prop:Sxi-isoinv}
Assume that, for a given $i\in I_p$, $K\times\{i\}\subset\bar{X}$ is an isolated invariant set for $\bar{f}^p$ in its isolating neighborhood $M\times\{i\}$. Then $K$ is an isolated invariant set for $f^p$, isolated by $M$.
\end{prop}
%------------------------------------------
\begin{prop}\label{prop:wip_ANR}
Let $f:\RR^d\to\RR^d$ be a discrete dynamical system. Set $\overline{\RR^d}:=\RR^d\times I_p$ and consider the dynamical system $\bar{f}$ on $\overline{\RR^d}$ given by (\ref{eq:def_fbar}).
Assume that $K:=\bigcup _{i=0} ^{p-1}\left(K_{\sigma_i}\times\{i\}\right)\subset \overline{\RR ^d}$ is an isolated invariant set with respect to $\bar{f}^p$, and
$M:=\bigcup _{i=0} ^{p-1}\left(M_{\sigma_i}\times\{i\}\right)$ is its isolating neighborhood. Then, there exists a weak index pair $Q$ for $\bar{f}^p$ and $K$ consisting of compact ANR's (for the definition of an ANR we refer to \cite{Bo}).
\end{prop}
\proof 
Fix an arbitrary $i\in I_p$. First note that $K_{\sigma_i}\times\{i\}=\Inv(M_{\sigma_i}\times\{i\},\bar{f}^p)$, as $\bar{f}^p$ maps $\RR^d\times\{i\}$ into itself. As a consequence, $M_{\sigma_i}\times\{i\}$ is an isolating neighborhood of $K_{\sigma_i}\times\{i\}$ with respect to $\bar{f}^p$. By Proposition \ref{prop:Sxi-isoinv}, $K_{\sigma_i}$ is an isolated invariant set with respect to $f^p$, and $M_{\sigma_i}$ is its isolating neighborhood. Using \cite[Lemma 5.1]{Szymczak-1996} we can take a polyhedral index pair $Q^{\sigma_i}$ for $f^p$ and $K_{\sigma_i}$. By \cite[Theorem 4.4]{Mr06}, $Q^{\sigma_i}$ is a weak index pair.
Then the pair $Q^{\sigma_i}\times\{i\}$ consists of compact ANR's, and constitutes a weak index pair for $\bar{f}^p$ and $K_{\sigma_i}\times\{i\}$. One can verify that the union $Q:=\bigcup_{i=0} ^{p-1} Q^{\sigma_i}\times\{i\}$ is a weak index pair with respect to $\bar{f}^p$ and $K$. Moreover, $Q_1$ and $Q_2$ are ANR's, as pairwise disjoint unions of ANR's.\qed
%------------------------------------------
\begin{thm}\label{thm:periodic_trajectory2}
Let $f:\RR^d\to\RR^d$ be a discrete dynamical system. Assume that $N=\bigcup_{i=1}^n N_i$, where $N_i$ are pairwise disjoint compact subsets of $N$, is an isolating neighborhood with respect to $f$, and $P$ is a weak index pair for $f$ in $N$. Let $p\in\NN$, let $\sigma:=(\sigma_0,\dots,\sigma_{p-1})\in I_n ^{\ZZ_p}$, and let endomorphism $I_\sigma$ of $\bigtimes _{i=0} ^{p-1}H^*(P^{\sigma_i})$ be given by
(\ref{eq:Isigma}). If
\begin{equation}\label{eq:Lef_Ip_neq0}
\Lambda (I_\sigma ^p)\neq 0
\end{equation}
then there exists a $p$-periodic point
$
x\in N_{\sigma_0}
$
for $f$ such that $f^{i+kp}(x)\in N_{\sigma _i}$, for $k\in\ZZ$.
\end{thm}
\proof
Consider the space $\overline{\RR^d}:=\RR^d\times I_p$, and the dynamical system $\bar{f}$ on $\overline{\RR^d}$, given by (\ref{eq:def_fbar}). By Proposition \ref{thm:periodic_trajectory1} we infer that $S_\sigma$ is an isolated invariant set with respect to $\bar{f}^p$. Thus, according to Proposition \ref{prop:wip_ANR}, we can take $Q$, a weak index pair for $\bar{f}^p$ and $S_\sigma$, consisting of compact ANR's. Then, by Proposition \ref{thm:periodic_trajectory1}, $\Lambda(I_{\bar{f}^p _{Q}})$ is well defined  and we have $\Lambda(I_{\bar{f}^p _{Q}})=\Lambda(I_{\sigma} ^p)$ which, along with 
(\ref{eq:Lef_Ip_neq0}), yields
$$
\Lambda(I_{\bar{f}^p _{Q}})\neq 0.
$$
Note that any weak index pair is a proper pair in the sense of \cite[Defnition 4]{Sr97}. Therefore, by \cite[Theorem 9]{Sr97}, there exists an $\bar{x}\in \cl(Q_1\setminus Q_2)$ such that $\bar{f}^p _{Q}(\bar{x})=\bar{x}$. Without loss of generality we may assume that $\bar{x}=(x,0)\in N_{\sigma_0}\times\{0\}$. Then, $x\in N_{\sigma_0}$ is a $p$-periodic point for $f$. Clearly $\{\bar{f}^k(\bar{x})\,|\,k\in\ZZ\}\subset S_\sigma$; hence $\{f^k(x)\,|\,k\in\ZZ\}\subset\Inv(\bigcup _{i=0} ^{p-1} N_{\sigma_i},f)$. Moreover, definition  (\ref{eq:def_fbar}) of $\bar{f}$ guaranties that the $p$-periodic trajectory of $f$ through $x$ passes through the components of $\Inv(\bigcup _{i=0} ^{p-1} N_{\sigma_i},f)$ in a proper order.
\qed\\
%%%%%%%%%%%%%%%%%%%%%%%%%%%%%%%%%%%%%%%%%%%%%%

We shall express the Lefschetz number of $I^p _\sigma$ in terms of the Lefschetz number of a composition of endomorphisms $g_i$ given by (\ref{eq:def_g}). Our goal is to prove the following theorem. 
%-------------------------------------------
\begin{thm}\label{thm:periodic_trajectory_g}
Let $f:\RR^d\to\RR^d$ be a discrete dynamical system. Assume that $N=\bigcup_{i=1}^n N_i$, where $N_i$ are pairwise disjoint compact subsets of $N$, is an isolating neighborhood with respect to $f$. Let $p\in\NN$, let $\sigma:=(\sigma_0,\dots,\sigma_{p-1})\in I_n ^{\ZZ_p}$, and let $P$ be a weak index pair for $f$ in $N$. Consider endomorphisms $g_i:H^*(P)\to H^*(P)$ given by (\ref{eq:def_g}). If
\begin{equation}\label{eq:Lef_comp_g_neq0}
\Lambda (g_{\sigma_0}\circ\dots\circ g_{\sigma _{p-1}})\neq 0
\end{equation}
then there exists a $p$-periodic point
$
x\in N_{\sigma_0}
$
for $f$ such that $f^{i+kp}(x)\in N_{\sigma _i}$, for $k\in\ZZ$.
\end{thm}
For the proof we need an auxiliary lemma.
%-----------------------------------------
\begin{lem}\label{lem:trIp_by_g}
Assume $I_\sigma$, $g_i$ and $h_i$ are given by (\ref{eq:Isigma}), (\ref{eq:def_g}) and (\ref{eq:def_h}), respectively. Then:
\begin{itemize}
\item[(i)] if $\Lambda (h_0\circ\dots\circ h_{p-1})$ is well defined then so is $\Lambda(I_\sigma ^p)$, and 
$$
\Lambda(I_\sigma ^p)=p\Lambda(h_{\sigma_0}\circ\dots\circ h_{\sigma_{p-1}}),
$$
\item[(ii)] if $\Lambda(g_{\sigma_0}\circ\dots\circ g_{\sigma_{p-1}})$ is well defined then so is $\Lambda(I_\sigma ^p)$ and we have
$$
\Lambda(I_\sigma ^p)=p\Lambda(g_{\sigma_0}\circ\dots\circ g_{\sigma_{p-1}}).
$$
\end{itemize}
\end{lem}
\proof
Note that $I_\sigma ^p$ and $h_{p-1}\circ\dots\circ h_{0}$ are endomorphisms of graded modules, however we consciously skip denoting the dimension in order to simplify the notation.
Observe that, by Lemma \ref{lem:Ip_h_g}(i) and the cyclic property of the trace, in each dimension we have the equality 
\begin{equation}\label{eq:trIp_by_h}
\tr(I_\sigma ^p)=p\tr(h_{\sigma_0}\circ\dots\circ h_{\sigma_{p-1}}).
\end{equation}
This completes the proof of (i).

For the proof of (ii) it suffices to verify that, in each dimension, we have
$$
\tr(I_\sigma ^p)=p\tr(g_{\sigma_0}\circ\dots\circ g_{\sigma_{p-1}}).
$$
Using Lemma \ref{lem:Ip_h_g}(ii), by the cyclic property of the trace, and (\ref{eq:g_comp=g}), we can write
$$
\begin{array}{rcl}
\tr(h_{\sigma_0}\circ\dots\circ h_{\sigma_{p-1}})&=&\tr(m_{\sigma_{p-1}}\circ\pi_{\sigma_{p-1}}\circ (g_{\sigma_0}\circ\dots\circ g_{\sigma_{p-1}})\circ \iota_{\sigma_{p-1}}\circ r_{\sigma_{p-1}})\\[1ex]
&=&\tr((g_{\sigma_0}\circ\dots\circ g_{\sigma_{p-1}})\circ \iota_{\sigma_{p-1}}\circ r_{\sigma_{p-1}}\circ m_{\sigma_{p-1}}\circ\pi_{\sigma_{p-1}})\\[1ex]
&=&\tr(g_{\sigma_0}\circ\dots\circ g_{\sigma_{p-1}}).
\end{array}
$$
Now, the assertion follows from (i).
\qed\\

{\bf Proof of Theorem \ref{thm:periodic_trajectory_g}.} The theorem follows from Theorem \ref{thm:periodic_trajectory2} and Lemma \ref{lem:trIp_by_g}.\qed

%%%%%%%%%%%%%%%%%%%%%%%%%%%%%%%%%%%%%%%%%%%
\section{Semiconjugacies to shift dynamics}
\label{sec:semiconjugacy}

Given a matrix $A\in\{0,1\}^{I_n\times I_n}$ we say that a partial map $s:\ZZ\pmap I_n$ is $A$-{\em admisssible}
if $A(s_i,s_{i+1})=1$ for any $i,i+1\in\dom s$.

Assume $V$ is a finite dimensional graded vector space over the field of rational numbers.
Let $V_i\subset V$ for $i\in I_n$ be subspaces of $V$ such that
$V=\oplus_{i=1}^nV_i$ is a direct sum decomposition of $V$ and let
\[
  p_i:V\ni x=(x_1,x_2,\ldots x_n)\longmapsto (0,0,\ldots 0,x_i,0\ldots,0)\in V_i.
\]
denote the canonical projections.

Consider a linear map $L:V\to V$. We define the {\em transition matrix} of $L$ with respect to the decomposition
$V=\oplus_{i=1}^nV_i$ as the matrix $A\in\{0,1\}^{I_n\times I_n}$ such that $A(i,j)=1$
if and only if $p_j\circ L\circ p_i\neq 0$. We say that $L$ is {\em Lefschetz-complete} if
\[
   \Lambda(L\circ p_{s_1}\circ L\circ p_{s_2}\circ \cdots \circ L\circ p_{s_k})\neq 0
\]
for any sequence $s:I_k\to I_p$ admissible with respect to the transition matrix of $L$.

Let $\Sigma_n:=\setof{s:\ZZ\to I_n}$ be the space of bi-infinite
sequences of elements in $I_n$ with product topology
and for a matrix $A\in\{0,1\}^{I_n\times I_n}$ let $\Sigma_A$ denote the subspace of
$A$-admissible sequnces. It is easy to see that the {\em shift map} $\sigma:\Sigma_n\to\Sigma_n$
defined by $\sigma(s)_{i}:=s_{i+1}$ is a homeomorphism and $\sigma(\Sigma_A)\subset\Sigma_A$.
Hence, $\sigma$ is a generator of a dynamical system on $\Sigma_A$.

\begin{thm}
\label{thm:decomposition}
Assume $N$ is an isolating neighborhood with respect to $f:\RR^d\to \RR^d$, and $P$ is a weak index pair for $f$ in $N$. Moreover, assume $N=\bigcup_{i=1}^n N_i$ where $N_i$ are pairwise disjoint compact subsets of $N$, and the index map $I_{f_P}:H^*(P)\to H^*(P)$ is Lefschetz-complete with respect to the decomposition $N=\bigcup_{i=1}^n N_i$. Then there exists a semiconjugacy $\rho$ between $S:=\Inv(\bigcup_{i=1}^n N_i,f)$ and the shift dynamics $\sigma$ on $\Sigma_A$, where $A$ is a transition matrix of $I_{f_P}$. Moreover,
for each periodic $s\in \Sigma_A$ there exists a periodic point of $f$ in $\rho ^{-1}(s)$.
\end{thm}
\proof
Fix an arbitrary $x\in S$. Since the sets $N_i$ are pairwise disjoint and $S=\Inv(\bigcup_{i=1}^n N_i,f)$, for each $k\in\ZZ$ there exists a unique $i\in I_n$ with $f^k(x)\in N_i$. By putting $\rho(x)_k:=i$ we define a continuous map $\rho:S\to 
\Sigma_{n}$. Note that, in fact, $\rho$ maps $S$ into $\Sigma_A$, as $\Sigma_{A}$ is the subspace of $\Sigma_{n}$ of all sequences admissible with respect to the transition matrix of $I_{f_P}$.     

We shall prove that $\rho$ is a surjection onto $\Sigma_{A}$. To this end let $s\in \Sigma_{A}$ be fixed. For an arbitrary $k\in\NN$ let $s^k$ denote the restriction of $s$ to the domain $\{-k,-k+1,\dots,k-1,k\}$. Since $I_{f_P}$ is  Lefschetz-complete,
we have 
$$
   \Lambda(I_{f_P}\circ p_{s_{-k}}\circ \cdots\circ I_{f_P}\circ p_{s_{0}}\circ \cdots \circ I_{f_P}\circ p_{s_k})\neq 0.
$$
By the cyclic property of the trace we obtain 
$$
   \Lambda(I_{f_P}\circ p_{s_{k}}\circ I_{f_P}\circ p_{s_{-k}}\circ \cdots\circ I_{f_P}\circ p_{s_{0}}\circ \cdots \circ I_{f_P}\circ p_{s_{k-1}})\neq 0,
$$
showing that $p_{s_{k}}\circ I_{f_P}\circ p_{s_{-k}}\neq 0$; hence, $(s_{-k},s_{k})$ is $A$-admissible. As a consequence, the periodic sequence $\tilde{s^k}:\ZZ\to I_n$, given by $\tilde{s^k}_m:=s^k _{(m+k)\mod (2k+1)-k}$ for $m\in\ZZ$, is $A$-admissible. By Theorem \ref{thm:periodic_trajectory_g}, there exists $x_k\in S$ such that $\rho(x_k)=\tilde{s^k}$. Since $k\in{\NN}$ was arbitrarily fixed, we have constructed a pair of sequences: $\{\tilde{s^k}\}\in\Sigma_{A} ^{\NN}$ convergent to $s$, and  $\{x_k\}\in S^{\NN}$, such that $\rho (x_k)=\tilde{s^k}$ for $k\in\NN$. By compactness of $S$, passing to a subsequence, if necessary, we may assume that $\{x_k\}$ converges to $x\in S$. Then, by the continuity of $\rho$ we have $\rho (x)=s$.  

The commutativity of the diagram
$$
\begin{tikzcd}
S\ar{d}[swap]{\rho}\ar{r}{f}&S\ar{d}{\rho}\\[3ex]
\Sigma_{A}\ar{r}[swap]{\sigma}&\Sigma_{A}\end{tikzcd}
$$
is easily readable. 

The above shows that $\rho$ constitutes a semiconjugacy from $f$ to the shift dynamics $\sigma$ on $\Sigma_{A}$.

The last statement of the theorem is a direct consequence of Theorem \ref{thm:periodic_trajectory_g}.\qed

%---------------------------------------------------------
Theorem \ref{thm:decomposition} has its counterpart in terms of endomorphisms $g_i$ given by (\ref{eq:def_g}).
\begin{thm}
\label{thm:decomposition_g}
Assume $N$ is an isolating neighborhood with respect to $f:\RR^d\to \RR^d$, and $P$ is a weak index pair for $f$ in $N$. Moreover, assume $N=\bigcup_{i=1}^n N_i$ where $N_i$ are pairwise disjoint compact subsets of $N$, and 
for each sequence $s:I_k\to I_p$ admissible with respect to the transition matrix $A$ of the index map $I_{f_P}:H^*(P)\to H^*(P)$ the composition $g_{s_1}\circ\cdots\circ g_{s_k}$ is not nilpotent. Then there exists a semiconjugacy $\rho$ between $S:=\Inv(\bigcup_{i=1}^n N_i,f)$ and the shift dynamics $\sigma$ on $\Sigma_A$.
\end{thm}
\proof
The proof runs along the lines of the proof of Theorem \ref{thm:decomposition}. Therefore, the details are left to the reader. However, it is worth to mention that now the admissibility of the periodic sequence $\tilde{s^k}:\ZZ\to I_n$ constructed in the proof of Theorem \ref{thm:decomposition} follows from the fact that  the composition $
 g_{s_{-k}}\circ \cdots\circ g_{s_{0}}\circ \cdots \circ g_{s_{k}}$ is not nilpotent. Moreover, the existence of the corresponding sequence $\{x_k\}\in S^{\NN}$ is guaranteed by Theorem \ref{thm:Wazewski_g}.
\qed
%%%%%%%%%%%%%%%%%%%%%%%%%%%%%%%%%%%%%%%%%%%

%%%%%%%%%%%%%%%%%%%%%%%%%%%%%%%%%%%%%%%%%%%%%%%%
\section{Proof of the main theorems}
\label{sec:henon}
\subsection{Proof of Theorem \ref{thm:henon3d}}
Clearly, $F$ is a cubical map. Its upper semicontinuity follows from \cite[Proposition 14.5]{G76}. Using elementary collapses  (cf. \cite{KMM}) we verify that $F$ has contractible values.

Using algorithms developed in \cite{Szymczak-1997}, a formula from \cite[Theorem 4.4]{B2017}, and techniques as in \cite{Pr19}, we find a cubical isolating block $N$ for $F$ consisting of five pairwise disjoint compact components $N_1,\ldots,N_5$, a cubical weak index pair $P$ in $N$, and index map $I_{F_P}$ (cf. Figure \ref{fig:henon3d}).
Direct computations show that $H^1(P_1,P_2) \cong \ZZ^5$ and $H^q(P_1,P_2) = 0$ for $q\neq 1$. More precisely, if we denote the generators of the cohomology group $H^1(P)$ by $\xi^1,\ldots,\xi^5$ and put $P^i:=P\cap N_i$, then we have
$$
H^1(P_1^i,P_2^i) =
\left\{\begin{array}{ll}
 \langle\xi^2\rangle& \mathrm{if}\ i=1,\\
 \langle\xi^5\rangle& \mathrm{if}\ i=2,\\
 \langle\xi^3\rangle& \mathrm{if}\ i=3,\\
 \langle\xi^1\rangle& \mathrm{if}\ i=4,\\
 \langle\xi^4\rangle& \mathrm{if}\ i=5.\\
\end{array}
\right.
$$
Moreover, using generators $\xi^1,\ldots,\xi^5$ as a basis, computations based on algorithms of \cite{MiMrPi05} provide the following matrix representation of the index map
$$
I_{F_{P}}^{1} = \left(
\begin{array}{rrrrr}
0 & 0 & -1 & 0 & 0 \\
-1 & 0 & 0 & 0 & 0 \\
0 & -1 & 0 & 0 & -1 \\
0 & 0 & 1 & 0 & 0 \\
0 & 0 & 0 & -1 & 0 \\
\end{array}
\right).
$$

Let $0<\varepsilon<\frac{1}{2}$ be fixed. By Theorem \ref{thm:index_maps_f_F_conj} we infer that there exists an $\varepsilon$-approximation of $F$, and each $\varepsilon$-approximation of $F$ shares with $F$ an isolating neighborhood and, up to a conjugacy, an index map.

Property (ii) is a straightforward consequence of Theorem \ref{thm:decomposition} under the assumption that $I_{F_{P}}$
is Lefschetz-complete. We verify this assumption by algorithmic computations. Details are presented in \cite{P2019}.

Finally, using the transition matrix $A$ we compute that the topological entropy of $f$ is greater than $\ln 1.2599$.
\qed

\subsection{Proof of Theorem \ref{thm:henon2d}}

The proof of this theorem is similar to the proof of previous theorem. We just note that the computations result in $H^1(P_1,P_2) \cong \ZZ^7$ and $H^q(P_1,P_2) = 0$ for $q\neq 1$. In particular, if $\xi^1,\ldots,\xi^7$ are generators of the cohomology group $H^1(P)$, then
$$
H^1(P_1^i,P_2^i) =
\left\{\begin{array}{ll}
 \langle\xi^2,\xi^6\rangle& \mathrm{if}\ i=1,\\
 \langle\xi^7\rangle& \mathrm{if}\ i=2,\\
 \langle\xi^4\rangle& \mathrm{if}\ i=3,\\
 \langle\xi^1\rangle& \mathrm{if}\ i=4,\\
 \langle\xi^5\rangle& \mathrm{if}\ i=5,\\
  \langle\xi^3\rangle& \mathrm{if}\ i=6,\\
\end{array}
\right.
$$
where $P^i:=P\cap N_i$, the components $N_i$ are pairwise disjoint and $N_i\subset N$ for each $i=1,\ldots,6$. We have also following matrix representation of the index map
$$
I_{F_{P}}^{1} = \left(
\begin{array}{rrrrrrr}
0 & 0 & 0 & 0 & 0 & 1 & 0 \\
0 & 0 & -1 & 0 & 0 & 0 & 0 \\
0 & 0 & 0 & -1 & 0 & 0 & -1 \\
0 & 0 & 0 & 0 & -1 & 0 & 0 \\
0 & 0 & 0 & 0 & 0 & 0 & -1 \\
0 & 0 & -1 & 0 & 0 & 0 & 0 \\
-1 & 0 & 0 & 0 & 0 & 0 & 0 \\
\end{array}
\right).
$$
The topological entropy of $\varepsilon$-approximation $f$ is greater than $\ln 1.151$.
\qed
%%%%%%%%%%%%%%%%%%%%%%%%%%%%%%%%%%%%%%%%%%%%%%%


\begin{thebibliography}{99}

\bibitem{B2017}
{\sc B.~Batko}.  Weak index pairs and the Conley index for discrete multivalued dynamical systems. Part II: properties of the Index, {\it SIAM J. Applied Dynamical Systems} {\bf 16} (2017), 1587--1617.

\bibitem{BM2016}
  {\sc B.~Batko and M.~Mrozek}.
  Weak index pairs and the Conley index for discrete multivalued dynamical systems, {\it SIAM J. Applied Dynamical Systems} {\bf 15} (2016), 1143--1162.

\bibitem{BMMP2018app}
  {\sc B.~Batko, K.~Mischaikow, M.~Mrozek and M.~Przybylski}.
  Conley index approach to sampled dynamics. Part II: applications. {\it preprint}

\bibitem{BEJM2019}
{\sc U.\ Bauer, H.\ Edelsbrunner, G. Jab\l{}o\'nski, and M.\ Mrozek}.
       Persistence in sampled dynamical systems faster, preprint 2017, arXiv:1709.04068 [math.AT]

\bibitem{Bo}
  {\sc K.~Borsuk}.
  Theory of retracs, PWN, Warszawa 1967.

\bibitem{C69}
{\sc A.~Cellina}.
Approximation of set-valued functions and fixed point theorems,
{\it Ann. Mat. Pura Appl.} {\bf 82} (1969), 17--24.

%\bibitem{C78}
%{\sc C. Conley}. Isolated invariant sets and the Morse index, {\it CBMS Lecture Notes}, {\bf 38} A.M.S. Providence, R.I. 1978.

\bibitem{D2003}
{\sc S.~Day}. A rigorous numerical method in infinite dimensions, {\it PhD diss.}, Georgia Institute of Technology (2003).

\bibitem{DFT2008}
{\sc S.~Day, R.~Frongillo, R.~Trevi{\~n}o}. Algorithms for Rigorous Entropy Bounds and Symbolic Dynamics, {\it SIAM J. Applied Dynamical Systems}, {\bf 7} (2008), 1477--1506.

\bibitem{dey2013}
{\sc T.\ Dey, F.\ Fan, Y.\ Wang}. Graph Induced Complex on Point Data, {\it Proceedings of the Twenty-ninth Annual Symposium on Computational Geometry}, SoCG '13 (2013), 107--116.

\bibitem{DJKKLM2019}
{\sc T.\ Dey, M.\ Juda, T.\ Kapela, J.\ Kubica, M. Lipi\'nski, and M.\ Mrozek}.
      Persistent Homology of Morse Decompositions in Combinatorial Dynamics.
      {\it SIAM Journal on Applied Dynamical Systems}, in print.

\bibitem{EJM2015}
{\sc H.\ Edelsbrunner, G. Jab\l{}o\'nski, M.\ Mrozek}.
         The Persistent Homology of a Self-map,
         {\it Foundations of Computational Mathematics}, {\bf 15}(2015), 1213–-1244.
         DOI: 10.1007/s10208-014-9223-y.

\bibitem{FMN2014}
{\sc S.~Ferry, K.~Mischaikow,V.~Nanda}. Reconstructing functions from random samples, {\it Journal of Computational Dynamics}, {\bf 1} (2014), 233-248.

\bibitem{Gr72}
{\sc A. Granas}. The Leray--Schauder index and the fixed point theory for arbitrary ANRs. {\it Bull. Soc. Math. France}, {\bf 100} (1972), 209--228.

\bibitem{G76}
{\sc L.~G\'orniewicz}.
Topological Fixed Point Theory of Multivalued Mappings, $2^{nd}$ ed.,
  {\it Topological Fixed Point Theory and Its Applications}
  {\bf 4}, Springer Verlag, The Netherlands, 2006.

\bibitem{GGK} {\sc L. G\'orniewicz, A. Granas, W. Kryszewski}, {\it On the homotopy method in the fixed point index theory of multi-valued mappings of compact absolute neighborhood retracts}, JMAA {\bf 161} (1991), 457--473.

\bibitem{HKMP2016}
{\sc S.~Harker, H.~Kokubu, K.~Mischaikow, and P.~Pilarczyk.}
Inducing a map on homology from a correspondence,
{\it Proc. AMS} {\bf 144} (2016), 1787--1801.

\bibitem{harker2018}
{\sc S.\ Harker, M.\ Kramar, R.\ Levanger, K. Mischaikow}. A Comparison Framework for Interleaved Persistence Modules, {\it arXiv:1801.06725} (2018).

\bibitem{KMM}
{\sc T.~Kaczynski, K.~Mischaikow, and M.~Mrozek.}
Computational Homology, {\it Applied Mathematical Sciences} 157, Springer-Verlag, 2004.

\bibitem{KM95}
{\sc T. Kaczynski, M. Mrozek}. Conley index for discrete multi-valued dynamical systems, {\it  Topol. \& Appl.} {\bf 65} (1995), 83-96.

\bibitem{KMV2005}
{\sc W.\ Kalies,  K.\ Mischaikow and R.\ VanderVorst}.
An Algorithmic Approach to Chain Recurrence,
{\it Found Comput Math}, {\bf 5} (2005)  409--449.
https://doi.org/10.1007/s10208-004-0163-9

\bibitem{Mi56}
{\sc E.~Michael}.
Continuous selections,
{\it Ann. Math.} {\bf 63} (1956), 361--382.

\bibitem{MiMrPi05}
{\sc K.~Mischaikow, M.~Mrozek, P.~Pilarczyk}. Graph Approach to the Computation of the Homology of Continuous Maps, {\it Foundations of Computational Mathematics}, {\bf 5.2} (2005), 199--229.

\bibitem{MiMrReSz99}
 {\sc K.~Mischaikow, M.~Mrozek, J.~Reiss, A.~Szymczak}. Construction of Symbolic
    Dynamics from Experimental Time Series, {\it Physical Review Letters}, {\bf 82} (1999), 1144--1147.

\bibitem{Mr06}
{\sc M.~Mrozek}. Index pairs algorithms, {\it Found. Comput. Math.}, {\bf 6} (2006), 457--493.

\bibitem{M90}
{\sc M.~Mrozek}. Leray functor and cohomological index for discrete dynamical systems, {\it TAMS.}, {\bf 318} (1990), 149--178.

\bibitem{Mr90}
 {\sc M.~Mrozek}. Open index pairs, the fixed point index and rationality of zeta functions, {\it Ergodic Theory Dynam. Systems}, {\bf 10} (1990), 555--564.

\bibitem{Mr91}
 {\sc M.~Mrozek}. Some remarks on Garay's conjecture, {\it Acta Mathematica Hungarica}, {\bf 57} (1991), 53--59.

  \bibitem{NSW2008}
 {\sc P.~Niyogi, S.~Smale, S.~Weinberger}. Finding the Homology of Submanifolds with High Confidence from Random Samples, {\it  Discrete {\&} Computational Geometry}, {\bf 39} (2008), 419--441.

 \bibitem{oudot2015}
 {\sc S.~Oudot, D.~Sheehy}. Zigzag Zoology: Rips Zigzags for Homology Inference, {\it Foundations of Computational Mathematics}, {\bf 15} (2015) 1151--1186.

 \bibitem{Pr19}
 {\sc M.~Przybylski}. Algorithmic computation of the Conley index for multivalued maps with no continuous selector in cubical spaces, {\it Schedae Informaticae},
 in print.

 \bibitem{P2019}
 {\sc M.~Przybylski}. PhD thesis, Jagiellonian University, Kraków 2019.

\bibitem{Sr97}
 {\sc R.~Srzednicki}. Generalized Lefschetz Theorem and fixed point index formula, {\it Topology and Its Applications}, {\bf 81} (1997), 207--224.

\bibitem{Szymczak-1995} {\sc A.\ Szymczak}. The Conley index for decompositions of isolated invariant sets, {\it Fund. Math.}, {\bf 148} (1995), 71--90.

\bibitem{Szymczak-1996}
{\sc A.\ Szymczak}. The Conley index and symbolic dynamics, {\it Topology}, {\bf 35} (1996), 287--299.

\bibitem{Szymczak-1997}
  {\sc A.\ Szymczak}. A combinatorial procedure for finding isolating neighborhoods and index pairs,
{\em Proceedings of the Royal Society of Edinburgh} {\bf 127A} (1997).

\bibitem{Szymczak-1999} {\sc A.\ Szymczak}. Index pairs: from dynamics to combinatorics and back,
{PhD Thesis}, Georgia Institute of Techmology, 1999.

\end{thebibliography}
\end{document}